\documentclass[a4paper, 10pt]{amsart}

\usepackage[latin1]{inputenc}
\usepackage{mathrsfs}
\usepackage{amssymb, amsmath}
\usepackage{amsthm}
\usepackage[english]{babel}

\usepackage[numeric, lite, initials, nobysame]{amsrefs}
\input xy
\xyoption{arrow} \xyoption{matrix} \xyoption{frame}%\OnlyOutlines
\xyoption{curve}

\newcommand{\Z}{{\mathbb Z}}

\newcommand{\K}{{\mathbb K}}\newcommand{\T}{{\mathbb T}}

\def\hara{\Ho\mathcal A_{R/A}}
\def\harf{\Ho\mathcal A_{R/F}}
\def\:{\colon}
\def\.{\cdot}
\def\<{\left\langle}
\def\>{\right\rangle}
\def\({\left(}
\def\){\right)}
\def\epsilon{\varepsilon}
\def\phi{\varphi}
\def\subset{\subseteq}

\def\leq{\leqslant}
\def\geq{\geqslant}

\def\lra{\longrightarrow}

\def\mapsto{\longmapsto}
\def\Mod{\mathsf{Mod}}

\def\dr{\mathscr{D}_R}

\def\smash{\wedge}

\def\xra{\xrightarrow}

\def\op{\text{op}}
\def\can{\text{can}}

\def\mr{\mathscr M_R}
\def\MU{{MU}}

\def\wt{\widetilde}
\def\wh{\widehat}
\def\xra{\xrightarrow}
\def\bda{b\mathscr D_A}
\def\bdb{b\mathscr D_B}
\def\bdf{b\mathscr D_F}
\def\SS{\mathbb S}
\def\FF{\mathcal F}

\def\GG{\mathcal G}

\def\FF{\mathcal F}
\def\II{\mathcal I}

\def\RR{\mathcal R}

\def\fenv{F^{ \text{e}}}

\def\dfenv{\mathscr{D}_{\fenv}}

\newtheorem{thm}{Theorem}
\newtheorem{lem}[thm]{Lemma}
\newtheorem{prop}[thm]{Proposition}\newtheorem{cor}[thm]{Corollary}

\numberwithin{equation}{section} \numberwithin{thm}{section}

\theoremstyle{remark}
\newtheorem{rem}[thm]{Remark}

\theoremstyle{definition}
\newtheorem{defn}[thm]{Definition}
\newtheorem{examples}[thm]{Examples}
\newtheorem*{conventions}{Notation and Conventions}
\newtheorem*{ack}{Acknowledgments}

\DeclareMathOperator{\Hom}{Hom} \DeclareMathOperator{\id}{id}
\DeclareMathOperator{\Ten}{T} \DeclareMathOperator{\Ho}{Ho}
\DeclareMathOperator{\End}{End} \DeclareMathOperator{\gr}{gr}
\DeclareMathOperator{\Ext}{Ext} \DeclareMathOperator{\Tor}{Tor}
\DeclareMathOperator{\im}{im} 
 \DeclareMathOperator{\coker}{coker}
\DeclareMathOperator{\holim}{holim}

\DeclareMathOperator{\Coext}{Coext}
\DeclareMathOperator{\Cohom}{Cohom} 
\DeclareMathOperator{\Sym}{Sym} 

 \DeclareMathOperator{\Der}{Der}
\DeclareMathOperator{\DDer}{\mathscr Der} 
 
 \DeclareMathOperator{\THH}{THH}
\DeclareMathOperator{\spec}{spec}

\newcommand{\ie}{i.e.}
\newcommand{\eg}{e.g.}

\title{Infinitesimal thickenings of Morava $K$-theories}
\author{Samuel W\"uthrich}
\date{23/03/07}
\address{EPFL SB IGAT, Batiment BCH, CH-1015 Lausanne, Switzerland}
\email{samuel.wuthrich@epfl.ch}


\begin{document}

\begin{abstract}
A.\@ Baker has constructed certain sequences of cohomology theories which interpolate
between the Johnson--Wilson and the Morava $K$-theories. We realize the representing
sequences of spectra as sequences of $\MU$-algebras. Starting with the fact that the
spectra representing the Johnson--Wilson and the Morava $K$-theories admit such
structures, we construct the sequences by inductively forming singular extensions. Our
methods apply to other pairs of $\MU$-algebras as well.
\end{abstract}

\maketitle

\section{Introduction}

The Johnson--Wilson theories and the Morava $K$-theories constitute
two families of multiplicative cohomology theories that play a
central role in the approach to stable homotopy theory known as the
chromatic point of view. It is an important fact that they can be
represented by $\SS$-algebras \cites{bj, rob}, meaning that the
homotopy ring structures on suitably chosen representing spectra
$E(n)$ and $K(n)$ (for a fixed prime $p$) can be rigidified in some
point-set category of spectra, like the ones described in \cite{ekmm}
or \cite{hss}, to strict ones. Moreover, there are maps of
$\SS$-algebras $\eta_n\: MU\to E(n)$ and $\rho_n\: E(n)\to K(n)$,
where $MU$ is a chosen commutative $\SS$-algebra representing complex
cobordism. It can be arranged that both $\eta_n$ and $\rho_n\eta_n$
are central maps, giving $E(n)$ and $K(n)$ the structure of
$MU$-algebras \cites{laz, bj}.

One application, which is the starting point for this paper, is Baker's construction
\cite{baker} of a sequence of $E(n)$-module spectra under $E(n)$ of the form
\begin{equation}\label{seqint}
E(n) \lra \cdots \lra E(n)/I_n^{s+1}\lra E(n)/I_n^s \lra \cdots \lra E(n)/I_n=K(n).
\end{equation}
Its homotopy limit is homotopy equivalent to $\widehat{E(n)}$, the
Bousfield localization of $E(n)$ with respect to $K(n)$, and the induced sequence of
homotopy groups is isomorphic to the canonical sequence of projections
\[
E(n)_* \lra \cdots \lra E(n)_*/I_n^{s+1}\lra E(n)_*/I_n^s \lra \cdots \lra E(n)_*/I_n
\cong K(n)_*,
\]
where $I_n$ denotes the kernel of the map induced by $\rho_n$ on coefficient rings.

The aim of this paper is to strengthen the bond between $E(n)$ and $K(n)$ that the sequence
\eqref{seqint} provides, by constructing it as a sequence of $\MU$-algebras. In fact, our
construction applies to a variety of other pairs of $\MU$-algebra spectra $T$ and $F$,
forming what we call a regular pair of $\MU$-algebras (Definition \ref{regularpair}).
This notion embraces a natural class of $\MU$-algebras $T$ and $F$ with commutative
coefficient rings which are related by an $\MU$-algebra map $T\to F$ that induces the
projection of $T_*$ onto a regular quotient $T_*/I\cong F_*$ on coefficients. Besides
$(E(n), K(n))$, examples are the pairs $(BP, P(n))$, $(\wh{E(n)}, K(n))$, $(E_n, K_n)$
and $({BP}\langle n\rangle, k(n))$ for an arbitrary prime $p$. The role of the ground
ring spectrum is not restricted to $\MU$. It can be played by any even commutative
$\SS$-algebra $R$, meaning that its homotopy groups are trivial in odd degrees.

\begin{thm}\label{theorem}
Let $(T, F)$ be a regular pair over an even commutative $\SS$-algebra $R$. There is a sequence of
$R$-algebra spectra under $T$ of the form
\begin{equation*}
T\lra \cdots \lra T/I^{s+1}\lra T/I^s \lra \cdots \lra T/I=F,
\end{equation*}
such that the induced sequence of coefficients is the canonical sequence of $T_*$-algebras
\[
T_* \lra \cdots \lra T_*/I^{s+1}\lra T_*/I^s \lra \cdots \lra T_*/I \cong F_*.
\]
\end{thm}

According to this result, the coefficients of $T/I^s$ are isomorphic to the ring of
functions on the scheme known as the infinitesimal thickening or neighbourhood of order
$s-1$ of the closed subscheme $\spec(T_*/I)\subseteq\spec(T_*)$ defined by $I$ (see \eg\
\cite{hartshorne}*{II.9} or \cite{illusie}*{1.3}). We therefore refer to the $R$-algebras
$T/I^s$ as infinitesimal thickenings of $F=T/I$.

The proof of Theorem \ref{theorem} relies on the results of \cite{sw}. There we construct
such ``topological $I$-adic towers'' as towers of $R$-module spectra. This is done in
such a way that $I^s/I^{s+1}$, the homotopy fibre of $T/I^{s+1}\to T/I^s$, is equivalent
to a wedge of suspensions of $F$. In the case where $T=R$ and $F$ is homotopy
commutative, topological $I$-adic towers were first constructed in \cite{bl}, using
different methods.

As an application of Theorem \ref{theorem}, we obtain a qualitative statement concerning
the relation of $T$ and $F$ if $T=R$. We briefly recall some definitions of
\cite{rognes}. A map of commutative $R$-algebras $A\to B$ is said to be symmetrically
\'etale if the canonical map $B\to \THH^A(B)$ to the topological Hochschild homology of $B$
relative to $A$ is an equivalence and $B$ is dualizable as
an $A$-module. For an algebra $D$ over a commutative $\SS$-algebra $C$, the structure map
$C\to D$ is defined to be symmetrically Henselian if each symmetrically \'etale map has
the unique lifting property up to contractible choice with respect to $C\to D$. Now let
$F$ be a regular quotient algebra of an even commutative $\SS$-algebra $R$, by which we
mean that $(R, F)$ is a regular pair of $R$-algebras. Consider the $R$-algebra $\wh
R=\holim_s R/I^s$. It may be interpreted as the Bousfield localization $L_F^R(R)$ of the
$R$-module spectrum $R$ with respect to $F$ \cites{bl, sw}. Therefore, it admits the
structure of a commutative $R$-algebra \cite{ekmm}*{Thm.\@ VIII.2.2}. Rognes
\cite{rognes}*{Prop.\@ 9.6.2} has proved the following statement for homotopy commutative $F$,
using the so-called external $I$-adic tower constructed in \cites{bl, laz}. The sequence
provided by Theorem \ref{theorem} leads to a new proof, which is valid for any regular quotient
algebra, not necessarily homotopy commutative.

\begin{cor}
Let $F$ be a regular quotient algebra of an even commutative $\SS$-algebra $R$. Then the
natural map $\wh R \to F$ is symmetrically Henselian.
\end{cor}

We construct the infinitesimal thickenings $T/I^s$ of $F=T/I$ inductively as singular
extensions. This method produces an $R$-algebra $B$ from ``strict derivations'' $A\to
A\vee M$, \ie\ maps of $R$-algebras over $A$ from a given $R$-algebra $A$ to the
square-zero extension $A\vee M$ of $A$ by some $A$-bimodule spectrum $M$. Homotopy
classes of such maps are possible to construct thanks to the fact that they correspond to
homotopy classes of $A$-bimodule maps from $D_A$, the fibre of the multiplication map
$A\smash_R A\to A$, to $M$ \cite{laz}. They are therefore strongly related to
$\THH_R^*(A, M)$, the topological Hochschild cohomology groups of $A$ with coefficients
in $M$.

As $I^s/I^{s+1}$ splits into a wedge of suspensions of $F$, the crucial calculation for
the inductive step of the construction is the identification of $\THH_R^*(T/I^s, F)$. We
will see that the problem can be reduced to the case where $T=R$. The situation is
fundamentally different for $s=1$ and $s>1$. For $s=1$, a certain universal coeffi\-cient
spectral sequence converging to $\THH^*_R(F, F)$ collapses for degree reasons, but the
extension problem is quite delicate. We show that $R/I^2$ can be constructed as a
singular extension of $F$, irrespective of the specific $R$-algebra structure of $F$.

For $s>1$, we first compute the endomorphism algebra $F^*_{R/I^s}(F)$ of $F$, viewed as a
left $R/I^s$-module spectrum, by means of a Bousfield--Kan spectral sequence. After that,
we construct a natural right coaction of $F^*_R(F)$, the coalgebra of $R$-linear
endomorphisms of $F$, on $F^*_{R/I^s}(F)$ and a natural map
\[
\THH_R^*(R/I^s, F) \lra P(F^*_{R/I^s}(F))
\]
into the coaction primitives of $F^*_{R/I^s}(F)$. We prove that it is an isomorphism,
using the fact that $F^*_{R/I^s}(F)$ is relatively injective over $F^*_R(F)$.

Our approach shows that sequences of infinitesimal thickenings exist, but does not
provide a way to classify all possible algebra structures on the $R/I^s$. For our
constructions, we have to make choices of certain strict derivations. It would be
desirable to know how different choices affect the obtained algebra structures. Essential
to pursuing this question is a thorough understanding of $\THH_R^*(R/I^s, F)$ for $s\geq
1$. As already mentioned, this problem is very different for $s=1$ and for $s>1$. Recent
results of Angeltveit \cite{anghochschild} concerning the former may help to classify
algebra structures on $R/I^2$ which are compatible with a given one on $F$. We compute
$\THH^*_R(R/I^s, F)$ for $s>1$ in this paper, but our identification is not natural and
not canonical. We hope to obtain a natural description in future work, by using
additional structure on $\THH^*_R(R/I^s, F)$.

The organization of the paper is as follows. In Section \ref{toptowers}, we recall the
results from \cite{sw} about the construction of topological $I$-adic towers. We add an
argument that shows that the construction is independent of the choice of a regular
sequence generating $I$. In Section \ref{singext}, we recall the definition and basic
facts about singular extensions of algebra spectra. In Section \ref{bockstein}, we
construct the first infinitesimal thickening of $F$. In Section \ref{thhcoh}, we relate
$\THH^*_R(G, F)$ and $F^*_G(F)$, for an $R$-algebra $G$ over a regular quotient algebra
$F$ of $R$, in the way indicated above. We characterize a situation which allows to
identify $F^*_G(F)$ and $\THH_R^*(G,F)$. In Section \ref{infthick}, we use these results
to construct the higher infinitesimal thickenings of $F$ inductively.

\begin{conventions}
We will use the results, notation and terminology from \cite{ekmm}. Throughout the paper,
we work over a fixed even commutative $\SS$-algebra $R$. We abbreviate $\smash_R$ by
$\smash$. Associated to an $R$-algebra $A$ are the strict and the derived categories
$\mathscr M_A$ and $\mathscr D_A$ of left $A$-module spectra, which have the same
objects. If we make a statement about an ``$A$-module'', we have in mind an $A$-module
spectrum viewed as an object in $\mathscr D_A$. A right $A$-module is a module over the
opposite $R$-algebra $A^\op$. If $B$ is another $R$-algebra, an $A$-$B$-bimodule is by
definition an $A\smash B^\op$-module. By an $A$-bimodule, we mean an $A$-$A$-bimodule. We
write $b\mathscr M_A$ and $b\mathscr D_A$ for the strict and derived categories of
$A$-bimodules. A free $A$-module is one equivalent to a wedge of suspensions of $A$. For
two $A$-modules $M$ and $N$, we write $\mathscr D_A^*(M, N)$ or $[M, N]^*_A$ for the
graded $R^*$-module of graded morphisms from $M$ to $N$ in $\mathscr D_A$. It coincides
with the homotopy groups $\pi_{-*}(F_A(M,N))$ of the internal function $R$-module
$F_A(M,N)$. Here, as well as throughout the paper, our convention is that we mean the
{\em derived} versions of functors like $F_A(-,-)$ whenever the arguments in question are
specified as $A$-modules. If we want to consider the strict functors, we clearly declare
the arguments as {\em strict} $A$-modules. If $N$ is an $A$-module and $M$ a right
$A$-module, we write $N^*_A(-)$ for the cohomological functor $\mathscr D_A^*(-, N)$ and
$M_*^A(-)$ for the homological functor $\pi_*(M\smash_A -)$, both defined on $\mathscr
D_A$. Recall that any object in the model categories $\mathscr{M}_A$ is fibrant.

An $R$-ring spectrum is a homotopy $R$-algebra, \ie\ a monoid in the homotopy category
$\dr$. The reduced homology and cohomology groups of an $R$-ring spectrum $B$ with
respect to an $R$-module $F$ are defined as
\begin{align*}
\wt F_*^R(B)  & = \coker(F_*\cong F_*^R(R) \xra{F_*^R(\eta)} F_*^R(B))
\\
\wt F^*_R(B)  & = \ker (F^*_R(B) \xra{F_R^*(\eta)} F^*_R(R)\cong F^*).
\end{align*}
Here, as well as later on, we use the letter $\eta$ to denote the unit of the $R$-algebra or
$R$-ring spectrum in question. The product map will be denoted by $\mu$.

We usually write $1$ for identity maps. For clarity, we sometimes denote it by the same
name as the object. An arrow labelled ``$\can$'' denotes a canonically given map, which
should be clear from the context. We follow the convention $X^*=X_{-*}$ for graded
objects. If $M_*$ and $N_*$ are modules over a graded commutative ring $\Lambda_*$, we
write $\Hom_{\Lambda_*}^*(M_*, N_*)$ and $\End^*_{\Lambda_*}(M_*)$ for the graded
$\Lambda_*$-linear homomorphisms and endomorphisms, respectively. The upper-case stars
refer to the cohomological grading of morphisms: An element $f\in\Hom^k(M_*, N_*)$ is a
homomorphism of degree $-k$, \ie\ $f_l\: M_l \to N_{l-k}$. The tensor, exterior and
symmetric algebras on $M_*$ are denoted by $\Ten^{*,*}_{\Lambda_*}(M_*)$,
$\Lambda^{*,*}_{\Lambda_*}(M_*)$ and $\Sym^{*,*}_{\Lambda_*}(M_*)$, respectively. One
grading comes from the construction of these algebras, the other one from $M_*$. We use
the same symbols with just one upper-case star to denote the corresponding algebras,
graded by the total gradation. We write $\Mod_{\Lambda_*}$ for the category of graded
modules over a graded ring $\Lambda_*$. An unadorned $\otimes$ means the tensor product
$\otimes_{R_*}$ over the coefficients of the $\SS$-algebra $R$.
\end{conventions}

\begin{ack}
The author would like to thank John Greenlees and Neil Strickland for many helpful
discussions, Alain Jeanneret for his support and encouragement and Srikanth Iyengar for
bringing the literature on Golod rings to his attention. He has been supported by a grant
of the Swiss National Science Foundation while this research was in progress. Part of it
was carried out during a very inspiring stay at the Mittag--Leffler institute. He would
like to express his gratitude to the organizers of the algebraic topology semester for
inviting him.
\end{ack}

\section{Topological $I$-adic towers}\label{toptowers}

Assume that $x=(x_1, \ldots, x_n)$ is a given regular sequence of elements in the
coefficient ring $R_*=\pi_*(R)$ of $R$, generating an ideal $J$. Then the homotopy groups
of the $R$-module $L_x$, defined as
\[
L_x = R/x_1 \smash R/x_2 \smash \cdots \smash R/x_n,
\]
are isomorphic to $R_*/J$. Here we use the notation $R/x_i$ to denote the cofibre of
$x_i$, viewed as a graded endomorphism of $R$ in the homotopy category $\dr$. If $x'$ is
another regular sequence generating the same ideal $J$, there is an equivalence $L_x \to
L_{x'}$ of $R$-modules \cite{ekmm}*{Cor.\@ V.2.10}. As $R$ does not have any endomorphisms of
odd degrees, the equivalence is canonical. So from now on, we will omit $x$ from the
notation and alternatively write $R/J$ for $L$.

Under the additional hypothesis that each element $x_i$ of a chosen sequence $x$ is a
non-zerodivisor, we can endow each $R/x_i$ with an $R$-ring structure
\cite{str}*{Prop.\@ 3.1}. Taken together, these products define one on $L$. Furthermore, it has
been observed in \cite{ang} that the argument used in \cite{rob} to prove that the Morava
$K$-theories admit $A_\infty$-structures can be generalized to show that in the present
situation $R/J$ can be promoted to an $R$-algebra.

\begin{defn}\label{regquotient}
We call $L=R/J$, endowed with a chosen $R$-ring structure, a {\em regular quotient ring}
of $R$. If in addition an $R$-algebra structure on $L$ is fixed, we refer to $L$ as a
{\em regular quotient algebra} of $R$.
\end{defn}

\begin{examples}
The spectrum $\wh{E(n)}$ representing completed Johnson--Wilson theory has a unique
commutative $\SS$-algebra structure \cite{br} for all $n$ and all primes $p$, so that the
Morava $K$-theory spectrum $K(n)$ can be realized as a regular quotient algebra of
$\wh{E(n)}$. Similarly, the Morava $E$-theory or Lubin--Tate spectrum $E_n$ is a
commutative $\SS$-algebra in a unique way for any $n$ and $p$ \cites{brob, gh}. The
version of Morava $K$-theory whose representing spectrum is commonly denoted by $K_n$
(see \eg\ \cite{ghmr}) can then be constructed as a regular quotient algebra of $E_n$.
\end{examples}

As mentioned in the introduction, the setup just described is more restrictive than
necessary and can be easily generalized to cover many more examples. Namely, assume that
in addition to $R$ and a regular sequence $x$ of non-zerodivisors generating an ideal
$J\lhd R_*$ as above, we are given an $R$-ring spectrum $T$ whose coefficient ring $T_*$ is
commutative. Moreover, we require that $x$ is regular on the $R_*$-module $T_*$ as well.
We write $I=J\cdot T_*\lhd T_*$. The $R$-module
\[
F = T/I = T\smash L = T\smash R/x_1\smash \cdots \smash R/x_n
\]
has homotopy groups $F_*\cong T_*\otimes R_*/J \cong T_*/I$. An $R$-ring structure on $L$
chosen as above and the given one on $T$ induce an $R$-ring structure on $F$.
Analogously, if $T$ is given as an $R$-algebra, any $R$-algebra structure on $L$ induces
one on $F$.

\begin{defn}\label{regularpair}
If $T$ is an $R$-ring spectrum satisfying the above conditions and an $R$-ring structure
on $L$ is fixed, we call $(T, F=T\smash L)$ a {\em regular pair of $R$-ring spectra}. If
moreover $T$ is an $R$-algebra and an $R$-algebra structure on $L$ is fixed, we say that
$(T, F)$ is a {\em regular pair of $R$-algebras}.
\end{defn}

\begin{rem}
In \cite{sw}, we considered regular pairs $(T, F=T\smash L)$ of $R$-ring spectra with the
additional property that the $R$-ring structure on $T$ is commutative and called such
triples $(R, T, F)$ regular. An inspection shows that the construction of topological
$I$-adic towers in \cite{sw} does not depend on the commutativity of $T$. The assumption
that $T_*$ is commutative is sufficient. In fact, $T$ enters into the construction only
at the end, as a certain Adams resolution of $R$ with respect to $L$ is smashed with $T$
over $R$. Commutativity of $T_*$ is used to identify the homotopy groups of the resulting
tower.
\end{rem}

\begin{examples}
The pairs $(E(n), K(n))$, $({BP}, P(n))$, $({BP}\langle n\rangle, k(n))$, $(\wh{E(n)},
K(n))$ and $(E_n, K_n)$ (see \cites{ravenel, ghmr}) are all regular pairs of
$\MU$-algebras, for all $n$ and $p$. For the first three pairs, this follows for instance
from \cite{ang}*{Thm.\@ 4.2} or is proved in \cite{bj}. For the last two pairs, it is a
consequence of the fact that $E(n)$ is an $\MU$-algebra, as the canonical map
$E(n)\to\wh{E(n)}$ is Bousfield localization with respect to $K(n)$ (see \eg\
\cite{sw}*{Cor.\@ 6.13}) and hence a map of $\MU$-algebras by \cite{ekmm}*{Thm.\@ VIII.2.1}, and
because the canonical map $\wh{E(n)}\to E_n$ is a map of commutative $\MU$-algebras
\cite{rognes}*{Thm.\@ 1.5}.
\end{examples}

In the following, $(T, F=T\smash L)$ denotes a regular pair of $R$-ring spectra, with
$F_*\cong T_*/I$ and $L_*\cong R_*/J$. The {\em Bockstein operations} $Q_1, \ldots,
Q_n\in L^*_R(L)$ are obtained by smashing
\[
\rho_i\beta_i\: R/x_i\lra \Sigma^{|x_i|+1} R/x_i
\]
with the identities on the other smash factors of $L=R/x_1\smash\cdots\smash R/x_n$. Here
$\beta_i$ and $\rho_i$ are taken from the cofibre sequences
\[
\Sigma^{|x_i|} R \xra{x_i} R \xra{\rho_i} R/x_i \xra{\beta_i} \Sigma^{|x_i|+1} R.
\]
The first part of the following result is \cite{str}*{Prop.\@ 4.15}, the second part is
proved as \cite{laz2}*{Lemma 2.6} or \cite{ang}*{Prop.\@ 4.1}.

\begin{prop}\label{homcohom}
There are isomorphisms of $L_*$-algebras
\[
L^*_R(L) \cong \Lambda_{L_*}(Q_1, \ldots, Q_n), \quad L_*^R(L^\op) \cong
\Lambda_{L_*}(a_1, \ldots, a_n),
\]
with $|a_i|=|x_i|+1$.
\end{prop}

By construction, the $Q_i$ depend on the choice of the sequence $x$. However, there is a
``coordinate-free'' description of $L^*_R(L)$, due to Strickland. We will use it to show
that the construction of $I$-adic towers does not depend on the choice of $x$. Recall
that a homotopy derivation $g\: L \to M$ to an $L$-bimodule spectrum $M$ is defined to be
a map of $R$-modules making the diagram
\[
\xymatrix{ L\smash L \ar[rr]^-{1\smash g \vee g\smash 1} \ar[d] &&
L\smash M \vee M\smash L\ar[d]
\\ L \ar[rr]^-g && M}
\]
commutative. The vertical maps are given by the product on $L$ and the biaction of $L$ on
$M$. We write $\DDer_R^k(L, M)$ for the set of derivations $L\to \Sigma^k M$. (We reserve
the notation $\Der_R^k(L, M)$ for strict derivations; see Section \ref{singext}.)
Strickland constructs a natural map $\DDer_R^*(L, L) \to \Hom_{R_*}^{*-1}(J/J^2, L_*)$
and shows that it is an isomorphism \cite{str}*{Cor.\@ 4.19}. The construction works for
any $L$-bimodule $M$. Strickland's result shows:

\begin{prop}\label{coordfree}
If $M$ is $L$-free, there is a natural isomorphism
\begin{equation}\label{hoderid}
\DDer_R^*(L, M) \cong \Hom_{R_*}^{*-1}(J/J^2, M_*).
\end{equation}
The endomorphism algebra $L^*_R(L)$ is the exterior algebra generated by $\DDer_R^*(L,
L)$.
\end{prop}

The construction of $I$-adic towers described in \cite{sw} can be slickened by using as
an auxiliary object a certain $R$-algebra $P$, whose homotopy ring realizes a polynomial
ring $R_*[y_1, \ldots, y_n]$ in generators $y_i$ of degree $|y_i|=|x_i|$. To define it,
we need the tensor $R$-algebra construction $\T_R(M)$ on an $R$-module $M$
\cite{ekmm}*{Thm.\@ VII.2.9}:
\[
\T_R(M) = R \vee M \vee (M\smash M) \vee \cdots.
\]
We let $P$ be the $R$-algebra
\[
P = \T_R(S_R^{|x_1|})\smash \cdots \smash \T_R(S_R^{|x_n|}),
\]
where $S_R^{k}$ is the cofibrant $k$-dimensional $R$-sphere. As a consequence of a
K\"unneth isomorphism, the homotopy ring $P_*$ is isomorphic to $R_*[y_1,\ldots, y_n]$.
Continuing the analogy between algebra and topology, we define $\gr_I(T)$ to be $F\smash
P$ and $\gr_J(R)$ to be $R/J\smash P=L\smash P$. The notation is legitimate, as the
homotopy groups of $\gr_I(T)$ realize the algebraic associated graded ring
\[
\gr_I(T_*)= T_*/I \oplus I/I^2\oplus I^2/I^3 \oplus \cdots
\]
and the ones of $\gr_J(R)$ realize $\gr_J(R_*)$. The $R$-algebra $P$ can be endowed with
a grading, by assigning the spheres $S_R^{|x_i|}$ degree one. By construction, it is a
graded free $R$-module spectrum. Furthermore, $\gr_I(T)$ and $\gr_J(R)$ are graded
bimodules over $P$. So it makes sense to speak of their components of homogeneous degree
\nolinebreak $s$. We abuse notation and denote them by $I^s/I^{s+1}$ and $J^s/J^{s+1}$,
respectively, as they have homotopy groups $I^s/I^{s+1}$ and $J^s/J^{s+1}$.

Now let $\delta_0=T\smash\delta'_0\: F\to \Sigma I/I^2$, where $\delta'_0\: L\to\Sigma J/J^2$
corresponds to the identity of $J/J^2$ under the isomorphism \eqref{hoderid} for $M=J/J^2$.
Consider
\[
F \xra{\delta_0} \Sigma I/I^2\xra{\can} \Sigma \gr_I(T).
\]
By extending scalars along the unit $R\to P$ of $P$, we obtain a map of graded
$P$-modules $\delta_*\: \gr_I(T) \to \Sigma \gr_I(T)$, with components
\[
\delta_s\: I^s/I^{s+1} \lra \Sigma I^{s+1}/I^{s+2}.
\]
These maps determine a topological $I$-adic tower, by means of relative homological
algebra in the triangulated category $\dr$ \cites{em, miller}. We briefly recall the
terminology and refer to \cite{sw} for more details and additional references.

A relative injective resolution of an $R$-module $M$ with respect to an $R$-ring spectrum
$E$ is a sequence of $R$-modules
\[
\ast \lra M \lra I_0 \lra I_1 \lra I_2 \lra \cdots
\]
with the following two properties. Firstly, each $I_k$ is a retract of an $R$-module of
the form $E\smash X_k$ for some $R$-module $X_k$. Secondly, the image of the sequence under the
functor $E\smash -$ is split, \ie\ equivalent to a sequence of the form
\[
\ast \lra Z_1 \lra Z_1 \vee Z_2 \lra Z_2 \vee Z_3 \lra \cdots,
\]
where the maps are the natural ones. The fundamental fact about relative injective
resolutions is that each such determines a unique (up to non-canonical equivalence)
diagram of the form (arrows with a circle denote maps of degree -1)
\begin{equation*}
\begin{array}{c}
\xymatrix@!=0.3cm{  M=M_0 \ar[dr]  && M_1 \ar[ll]|-\circ\ar[dr] && M_2 \ar[ll]|-\circ\ar[dr]
&& M_3 \ar[ll]|-\circ&\ar@{}[d]|{.\,.\,.}
\\
& I_0 \ar[ur] \ar[ur] && I_1 \ar[ur] \ar[ur] && I_2 \ar[ur] & & }
\end{array}
\end{equation*}
which is composed of cofibre sequences
\[
M_i \lra I_i \lra M_{i+1} \lra \Sigma M_i.
\]
Let $i\:I^{s+1} \to I^{s}$ and $p\: I^s\to I^s/I^{s+1}$ denote the canonical inclusions and
projections respectively.

\begin{thm}[\cite{sw}*{Thm.\@ 6.6}]\label{iadic}
The sequence of $R$-modules
\begin{equation}\label{resseq}
\ast \lra T \xra{\eta} F=T/I \xra{\delta_0} \Sigma I/I^2 \xra{\delta_1} \Sigma^2 I^2/I^3 \lra
\cdots
\end{equation}
is a relative injective resolution of\/ $T$ with respect to $F$. The homotopy groups of
its associated Adams tower
\begin{equation}\label{iadictower}\begin{array}{c}
\xymatrix@!C=.8cm{ \cdots\ \ar[r]^{\iota} & I^3 \ar[r]^{\iota}\ar[d]_{\pi} & I^2
\ar[r]^{\iota}\ar[d]_{\pi} & I \ar[r]^{\iota}\ar[d]_{\pi} & T\ar[d]_{\pi}
\\
& I^3/I^4\ar[lu]^{\epsilon^3}|-\circ & I^2/I^3\ar[lu]^{\epsilon^2}|-\circ &
I/I^2\ar[lu]^{\epsilon^1}|-\circ & T/I \ar[lu]^{\epsilon^0}|-\circ}
\end{array}
\end{equation}
realize the $I$-adic filtration of $T_*$
\begin{equation*}
\xymatrix@!C=.8cm{ \cdots\ \ar[r]^{i} & I^3 \ar[r]^{i}\ar[d]_{p} & I^2
\ar[r]^{i}\ar[d]_{p} & I \ar[r]^{i}\ar[d]_{p} & T_*. \ar[d]_{p}
\\
& I^3/I^4 & I^2/I^3 & I/I^2\ & T_*/I }
\end{equation*}
\end{thm}

By ``reversing'' the tower in the statement of the theorem, using the octahedral axiom
\cite{hps}*{Def.\@ A.1.1}, we obtain a tower
\begin{equation}\label{quotiadic}\begin{array}{c}
\xymatrix{ T \quad \ar[r] & \quad\ldots\quad \ar[r]^-\rho & T/I^3 \ar[r]^-\rho & T/I^2
\ar[dl]^-{\theta^2} \ar[r]^-\rho & T/I \ar[dl]^-{\theta^1}
\\
&& I^2/I^3 \ar[u]^-\nu & I/I^2 \ar[u]^-\nu}
\end{array}
\end{equation}
whose homotopy groups are isomorphic to what the notation suggests.

\begin{prop}[\cite{sw}*{Thm.\@ 6.6}]\label{ses}
The long exact sequences of $F^*_R(-)$-cohomology groups derived from \eqref{quotiadic}
induce short exact sequence of free $F^*$-modules
\[
0 \lra \wt{F}_R^*(T/I^{s+1}) \lra F_R^*(I^s/I^{s+1}) \lra \wt{F}^{*+1}_R(T/I^s) \lra 0.
\]
\end{prop}

\begin{rem}\label{remt}
If $T$ is an $R$-algebra, we have
\[
F^*_T(F) \cong F^*_R(L) \cong\Lambda_{F^*}(Q_1, \ldots, Q_n).
\]
The proof of \cite{sw}*{Thm.\@ 6.6} shows that there are short exact sequences of free
$F^*$-modules
\begin{equation}\label{sest}
0 \lra \wt{F}_T^*(T/I^{s+1}) \lra F_T^*(I^s/I^{s+1}) \lra \wt{F}^{*+1}_T(T/I^s) \lra 0
\end{equation}
and that the sequences of Proposition \ref{ses} are isomorphic to the image of
\eqref{sest} under $F^*_R(T)\otimes-$.
\end{rem}

\begin{rem}
The constructions described in this section can be generalized to cover quotients by ideals
generated by infinite sequences. See \cite{sw}.
\end{rem}

\begin{rem}
An alternative construction of topological $I$-adic towers for the case where $T=R$ and
where $F$ is homotopy commutative is given in \cite{bl}.
\end{rem}

\section{Singular extensions of $\SS$-algebras}\label{singext}

Let $A$ be an $R$-algebra, not necessarily commutative. An $R$-algebra over $A$ is an
$R$-algebra $B$ together with a map $\pi\: B\to A$ of $R$-algebras. A map of $R$-algebras
over $A$ is a map of $R$-algebras which respects the structure maps. We write $\mathcal
A_{R/A}$ for the resulting category of $R$-algebras over $A$. The model structure on the
category of $R$-algebras gives rise to one on $\mathcal A_{R/A}$ \cite{hovey}*{Prop.\@ 1.1}.
We denote its homotopy category by $\hara$.

Associating to an $A$-bimodule $M$ the square zero extension $A\vee M$ defines a functor
$b\mathscr M_A\to\mathcal A_{R/A}$. Note that the inclusion $i\: A\to A\vee M$ is map in
$\mathcal A_{R/A}$. If $B$ is another $R$-algebra over $A$ and $f\: B\to A\vee M$ a map
in $\hara$, we can form the homotopy pullback of $f$ and $i$ in $\mathcal A_{R/A}$:
\begin{equation}\label{pullback}
\begin{array}{c}
\xymatrix{ C \ar[r]\ar[d] &  B\ar[d]^f \\
A\ar[r]^-i &   A\vee M}
\end{array}
\end{equation}

\begin{defn}
The $R$-algebra $C$ together with the map $C\to B$ from \eqref{pullback} is the {\em
singular extension} associated to $f$.
\end{defn}

By neglect of structure, \eqref{pullback} defines a homotopy pullback diagram in
$b\mathscr M_C$ \cite{dhks}*{19.6}. Applying to this the forgetful functors $b\mathscr
M_C\to\mathscr M_C\to\mr$ gives homotopy pullback diagrams in $\mathscr M_C$ and $\mr$
respectively. As $b\mathscr D_C$, $\mathscr D_C$ and $\dr$ are all additive, $C$ is
equivalent to the homotopy fibre of the composition $\tilde f\: B\xra{f}B\vee M\xra{\can}
M$ in those categories. Thus, we obtain a cofibre sequence of $C$-bimodule spectra
\[
C \lra B \xra{\ \ \tilde f\ \ } M\lra \Sigma M,
\]
which maps to cofibre sequences in $\mathscr D_C$ and $\dr$ under the forgetful functors.

The $B$-bimodule of associative differentials $D_B$ is defined to be the homotopy fibre
of the multiplication $B\smash B\to B$. (We use the notation $D_B$ rather than $\Omega_B$
from \cite{laz}, to distinguish the notion from the module of commutative differentials
of a commutative $R$-algebra, defined in \cite{bas}.) Extending scalars along $B\to A$
yields the $A$-bimodule $D_B^A=A\smash_B \wt D_B \smash_B A$, where $\wt D_B$ denotes a
cofibrant replacement of $D_B$ as a $B$-bimodule. According to our convention that smash
products are to be taken in the derived sense, we would write this simply as $A\smash_B
D_B \smash_B A$. We deviate from our convention here to avoid potential confusion. For an
$A$-bimodule $M$, the $R$-module of derivations from $B$ to $M$ is defined as
\[
\Der_R(B, M) = F_{A\smash A^\op}(D_B^A, M).
\]
We refer to the elements of
\[
\Der_R^{-k}(B,M) = \pi_k(\Der_R(B, M)) = \bda^{-k}(D_B^A, M)
\]
as strict derivations from $B$ to $M$ of degree $k$. This terminology is justified by the
following result of Lazarev. The grading on the right hand side is obtained by considering maps into
the square zero extensions $A\vee\Sigma^k M$, for $k\in\Z$.

\begin{thm}[\cite{laz}*{Thm.\@ 2.2}]\label{deradj}
There is an adjunction
\[
\Der_R^*(B, M) \cong \hara^*(B, A\vee M).
\]
\end{thm}

Recall that the topological Hochschild homology and cohomology groups of $B$ with
coefficients in a $B$-bimodule $N$ are defined as
\begin{equation*}
\THH_*^R(B, N) = \pi_*(\wt B\smash_{B\smash B^\op} N), \quad \THH^*_R(B, N) = \bdb^*(\wt
B, N),
\end{equation*}
respectively, where $\wt B$ denotes a cofibrant replacement of $B$ as a $B$-bimodule.
Regarding an $A$-bimodule $M$ as a $B$-bimodule via restriction along $B\to A$, we obtain
a long exact sequence of $R^*$-modules
\begin{equation}\label{les}
\cdots\to \THH^*_R(B, M) \to  M^* \to \Der^*_R(B, M) \to
\THH^{*+1}_R(B, M)\to \cdots.
\end{equation}
The ``universal derivation'' is defined to be the map of $R$-modules% $B$-bimodules
\[
d_B\: B\xra{\eta_B} A\vee D_B^A \xra{\can} D_B^A,
\]
where $\eta_B$ is the unit of the adjunction. Via the forgetful functor $U\:\bda\to\dr$,
$d_B$ induces a ``forgetful map''
\[
V_B^*\: \Der_R^*(B, M) \xra{U} \dr^*(D_B^A, M)\xra{d_B^*} \dr^*(B, M) = M^*_R(B).
\]
If $\tilde f\: B\to M$ is a given map of $R$-modules and $\phi\in\Der_R^0(B,M)$ is chosen
so that $V_B^*(\phi)=f$, we say that ``$\phi$ realizes $f$ as a strict derivation''.

By neglect of structure, the $R$-algebra map $f\: B\to A\vee M$ over $A$ corresponding to
a strict derivation $\phi\: D^A_B\to M$ defines a map of $R$-ring spectra over $A$. By
definition, this is nothing but a homotopy derivation in the sense of Section
\ref{toptowers}. This observation easily implies:

\begin{lem}\label{derivation}
The image of the forgetful map $V_B^*\: \Der_R^*(B, M) \to M^*_R(B)$ is contained
in $\DDer^*_R(B, M)$. In particular, we always have $\im(V_B^*)\subseteq \wt M^*_R(B)$.
\end{lem}

As pointed out in \cite{laz}, there is an analogue of Theorem \ref{deradj} for
derivations over a not necessarily commutative $R$-algebra $T$. Instead of $\mathcal
A_{R/A}$, we consider the category $\mathcal A_{T/A}$ of $T$-algebras over $A$. Its
objects are $R$-algebras $B$ with $R$-algebra maps $T\to B$ and $B\to A$. The morphisms
are the $R$-algebra maps compatible with these structure maps. For similar reasons as
before, $\mathcal A_{T/A}$ is a model category, whose homotopy category we denote by $\Ho
\mathcal A_{T/A}$. The product map $B\smash B\to B$ of a $T$-algebra $B$ over $A$ factors
through $B\smash_T B$. This makes $B$ a monoid in the monoidal category of $T$-bimodules,
with unit given by the structure map $T\to B$. The $B$-bimodule $D_{B/T}$ is defined as
the homotopy fibre of $B\smash_T B \to B$ and $D_{B/T}^A$ as $A\smash_B \wt D_{B/T}
\smash_B A$, with $\wt D_{B/T}$ a cofibrant replacement of the $B$-bimodule $D_{B/T}$.
The analogue to Theorem \ref{deradj} is then \cite{laz}*{Rem.\@ 2.4}
\[
\bda^*(D^A_{B/T}, M) \cong \Ho\mathcal A_{T/A}^*(B, A\vee M).
\]
We call the elements of $\bda^*(D^A_{B/T}, M)$ strict derivations from $B$ to $M$ over
$T$.

\section{Constructing the first order infinitesimal thickening}\label{bockstein}

From now on, we assume that $(T, F=T\smash L)$ is a regular pair of $R$-algebras, with
$F_*\cong T_*/I$, where $I=(x_1, \ldots, x_n)\cdot T_*$ for some regular sequence $(x_1,
\ldots, x_n)$ in $R_*$. The first aim of this section is to show that the homotopy
derivation $\delta_0\: F\to\Sigma I/I^2$ defined in Section \ref{toptowers} can always be
realized as a strict derivation, independently of the chosen $R$-algebra structure on
$F$. This has been shown by Lazarev for $R=T=\MU$ and $F=H\Z/p$ \cite{laz}*{Thm.\@ 10.2}.
The proof for the general case is analogous. The argument given in \cite{laz} is rather
brief, however, and it has taken the author some time to work out a complete proof. We
therefore include it here, for the convenience of the reader.

Let us first assume that $R=T$, so $F=L$. As a consequence of Proposition
\ref{homcohom}, the $E_2$-term of the universal coefficient spectral sequence
\cite{ekmm}*{IV.4}
\begin{equation}\label{thhseq}
E_2^{*,*} = \Ext_{F_{-*}^R(F^{\op})}^{*,*}(F^*, F^*) \Longrightarrow \THH^*_R(F, F)
\end{equation}
is isomorphic to a polynomial algebra $F_*[z_1, \ldots, z_n]$ in elements $z_i$ of
bidegree $(1, -|x_i|-1)$. It will be more convenient to have a description
which does not depend of the choice of a basis of $F_*^R(F^\op)$. Let us denote the
augmentation ideal $\ker(\mu_*\: F_*^R(F^\op) \to F_*)\cong (D_F)_*$ by $D$ and consider the
short exact sequence of $F_*^R(F^\op)$-modules
\[
\mathcal E\: \quad 0\lra D/D^2 \lra F_*^R(F^\op)/D^2 \lra F_* \lra 0.
\]
Note that the we may naturally identify
\[
(D/D^2)^\vee := \Hom_{F_*}(D/D^2, F_*) \cong \Hom_{F_*^R(F^\op)}(D/D^2, F_*).
\]
Thus, the connecting homomorphism associated to $\mathcal E$ defines a map
\begin{equation}\label{extmap}
(D/D^2)^\vee \lra \Ext_{F_*^R(F^\op)}^{1,*}(F_*, F_*).
\end{equation}
It induces an algebra isomorphism
\begin{equation}\label{isosym}
\Sym^{*,*}_{F_*}((D/D^2)^\vee) \cong \Ext_{F_{*}^R(F^\op)}^{*,*}(F_*, F_*)
\end{equation}
of the sort we were looking for.
Now it is clear that the spectral sequence \eqref{thhseq} collapses for degree
reasons. Hence, we find that the associated graded of $\THH^*_R(F,F)$, with respect to the
spectral sequence filtration, is given by
\[
\gr^*(\THH^*_R(F,F)) \cong \Sym^{*}_{F^*}((D/D^2)^\vee).
\]
Using the identification
\begin{equation}\label{derid}
\Der^*_R(F, \Sigma I/I^2) \cong \ker(\THH^*_R(F, \Sigma I/I^2)\lra
(\Sigma I/I^2)^*),
\end{equation}
and noting that $D/D^2\cong(\Sigma I/I^2)_*$, we then have
\begin{equation}\label{derfiso}
\gr^*(\Der^*_R(F, \Sigma I/I^2)) \cong \wt\Sym^*_{F^*}(\End_{F^*}(D/D^2)),
\end{equation}
where $\wt\Sym^*_{F^*}(\End_{F^*}(D/D^2))$ is the augmentation ideal. The solution of the
extension problem depends on the $R$-algebra structure chosen on $F$
\cite{anghochschild}. For our purposes, it suffices to know $\Der^*_R(F, \Sigma I/I^2)$
up to associated graded.

\begin{lem}\label{cofhom}
The sequence of homotopy groups of the cofibre sequence associated to the universal
derivation $d_F\: F\to D_F$ is of the form
\begin{equation}\label{sesdf}
0 \lra D_F^*\cong D/D^2\oplus D^2 \xra{i\oplus 1} R^{*+1}/I^2 \oplus D^2 \lra F^{*+1}
\lra 0,
\end{equation}
where $i$ is the canonical inclusion $I/I^2\to R^*/I^2$, with degree shifted by one.
\end{lem}

\begin{proof}
By definition, $d_F$ is the image of the identity on $D_F$ under the forgetful map
\[
\bdf(D_F, D_F) \cong \Der_R^*(F, D_F) \xra{V^*_F} \dr^*(F, D_F).
\]
It is therefore a homotopy derivation, by Lemma \ref{derivation}. Now consider the
universal coefficient spectral sequence
\begin{equation}\label{ucssdf}
E_2^{*,*} = \Ext_{R^*}^{*,*}(F^*, D_F^*) \ \Longrightarrow \ \dr^*(F, D_F).
\end{equation}
It was noted above that $D_F$ is equivalent to a free $F$-module spectrum. Hence, as the
analogous spectral sequence converging to $F^*_R(F)$ collapses, so does \eqref{ucssdf}.
It now suffices to show that the element representing $d_F$ lies in the first
$\Ext$-group and corresponds to the short exact sequence \eqref{sesdf} (see \eg\
\cite{sw}*{Prop.\@ 6.3}).

As $F_*^R(F)$ is $F_*$-free, the composition
\begin{equation}\label{nattrdf}
\dr^*(F, D_F) \lra \Hom^*_{F_*}(F_*^R(F), F_*^R(D_F)) \lra \Hom^*_{F_*}(F_*^R(F),
(D_F)_*)
\end{equation}
induced by applying first $F_*^R(-)$ and then using the left action of $F$ on $D_F$, is
an isomorphism (see \eg\ \cite{sw}*{Lemma 6.2}). Now consider the following diagram,
where the top map is \eqref{nattrdf}, the left map is the natural inclusion, the right
map is induced by the projection $F_*^R(F) \to D/D^2$ and the lower map is the
isomorphism \eqref{hoderid} from Proposition \ref{coordfree}:
\begin{equation}\begin{array}{c}\label{diader}
\xymatrix{\dr^*(F, D_F) \ar[r]^-\cong & \Hom^*_{F_*}(F_*^R(F), (D_F)_*)
\\
\DDer^*_R(F, D_F)\ar[u]\ar[r]^-\cong & \Hom^*_{F_*}(D/D^2, (D_F)_*).\ar[u]}
\end{array}
\end{equation}
We leave it to the reader to verify that the diagram commutes. Alternatively, he may
check that, for an $F$-free $F$-bimodule $M$, the restriction of the natural isomorphism
\cite{sw}*{Lemma 6.2}
\[
\dr^*(F, M) \cong \Hom_{F_*}^*(F_*^R(F), M_*)
\]
to $\DDer^*_R(F, M)$ factors through $\Hom^*_{F_*}(D/D^2, M_*)$, and that the induced map
is an isomorphism. This observation leads to another construction of the isomorphism
\eqref{hoderid}.

We claim that the homomorphism $F_*^R(F)\to (D_F)_*$ corresponding to $d_F$, viewed as an
$R$-module map, is given by the projection
\begin{equation*}%\label{proj}
F_*^R(F)\cong (D_F)_* \oplus F_* \lra (D_F)_*.
\end{equation*}
To show this, we use the fact \cite{laz}*{Thm.\@ 4.1} that the composition $F\xra{d_F}
D_F \xra{\can} F\smash F$ coincides with $\eta\smash F-F\smash\eta$. This map corresponds
under the isomorphism
\[
\dr^*(F, F\smash F) \cong \Hom^*_{F_*}(F_*^R(F), F_*^R(F))
\]
to the homomorphism $x\mapsto x-\mu_*(x)\cdot 1$. By naturality, this proves our claim.
As diagram \eqref{diader} commutes, it follows that the homotopy derivation $d_F$
corresponds to the inclusion $D/D^2\to (D_F)_*$. Now the short exact sequence
\[
0\lra D/D^2 \lra R^{*+1}/I^2 \lra F^{*+1} \lra 0
\]
gives rise to a natural isomorphism
\[
\gamma\: \Hom^*_{F_*}(D/D^2, (D_F)_*) \cong \Ext_{R_*}^{1,*}(F_{*-1}, (D_F)_*).
\]
It is not difficult to see that the induced isomorphism
\[
\DDer^*_R(F, D_F)\cong \Ext_{R_*}^{1,*}(F_{*-1}, (D_F)_*)
\]
maps a given derivation to the element which represents it in the spectral sequence
\eqref{ucssdf}. As $\gamma$ sends the inclusion $D/D^2\to (D_F)_*$ to \eqref{sesdf}, the
lemma is proved.
\end{proof}

\begin{prop}\label{strictbs}
Any element in $\Der^*_R(F, \Sigma I/I^2)$ lying in the coset corresponding under
\eqref{derfiso} to the identity of $D/D^2$ realizes $\delta_0\: F\to\Sigma I/I^2$ as a
strict derivation.
\end{prop}

\begin{proof}
We have to check that the forgetful map $V_F^*$, defined as the composition
\[
V_F^*\: \Der^*_R(F, \Sigma I/I^2) \xra{U} \dr^*(D_F, \Sigma I/I^2) \xra{d_F^*}
\dr^*(F, \Sigma I/I^2),
\]
maps an element in the stated coset to $\delta_0$. To that end, we consider the morphisms
induced by $U$ and $d_F^*$ on the universal coefficient spectral sequences
\begin{equation*}%\array{c}\label{specseq}
\xymatrix@C=0cm@R=0.5cm{ (E_2^{*,*})_1 &=&
\Ext_{F_{-*}^R(F^\op)}^{*,*}(D_F^*, (\Sigma I/I^2)^*)\ar[d]^-{U^*} & \Longrightarrow &
\Der_R^*(F,\Sigma I/I^2)\ar[d]^-{U}
\\
(E_2^{*,*})_2 &=& \Ext_{R^*}^{*,*}(D_F^*, (\Sigma I/I^2)^*)\ar[d]^-{d_F^*} &
\Longrightarrow & \dr^*(D_F,\Sigma I/I^2)\ar[d]^-{d_F^*}
\\
(E_2^{*,*})_3 &=& \Ext_{R^*}^{*,*}(F^*, (\Sigma I/I^2)^*) & \Longrightarrow &
\dr^*(F,\Sigma I/I^2).}
%\end{array}
\end{equation*}
All three spectral sequences collapse: We have seen this for the first one above; for the
second and the third one this is a consequence of \cite{sw}*{Prop.\@ 6.5}. Here we use
the fact that $D_F$, as an $F$-module spectrum with $F^*$-free homotopy groups, is
equivalent to a free $F$-module spectrum.

By Lemma \ref{cofhom} below, $d_F^*\: F^*\to D_F^*$ is trivial. Hence the induced morphism
on the associated graded
\[
\gr(d_F^*)\: \gr^*(\dr^*(D_F, \Sigma I/I^2)) \lra \gr^*(\dr^*(F, \Sigma I/I^2))
\]
is zero and so gives rise to a map
\[
\wt\gr^*(d_F^*)\: \gr^*(\dr^*(D_F, \Sigma I/I^2)) \lra \gr^{*+1}(\dr^*(F, \Sigma I/I^2)).
\]
We now check that the composite
\[
\wt\gr^*(d_F^*)\circ\gr^*(U)\: \gr^*(\Der^*_R(F, \Sigma I/I^2)) \lra \gr^{*+1}(\dr^*(F,
\Sigma I/I^2))
\]
maps the coset $x$ corresponding to the identity of $D/D^2$ to the coset representing
$\delta_0$. The element of $\Ext_{F_{-*}^R(F^\op)}^*(D_F^*, D/D^2)$ corresponding under
\[
\Ext_{F_{-*}^R(F^\op)}^*(D_F^*, D/D^2) \cong \gr^*(\Der^*_R(F, \Sigma I/I^2))
\]
to $x$ is the projection $q\: D_F^* \cong D\to D/D^2$. By the Geometric Boundary Theorem
(see \cite{sw}*{Prop.\@ 6.4} for instance) and Lemma \ref{cofhom} below,
$\wt\gr^*(d_F^*)$ is represented by the connecting homomorphism of $\Ext_{R^*}(-, F^*)$
associated to the short exact sequence \eqref{sesdf}. It maps $q$ to the short exact
sequence
\[
0\lra D/D^2 \lra R^{*+1}/I^2 \lra F^{*+1} \lra 0,
\]
which, by \cite{sw}*{Prop.\@ 6.5}, represents $\delta_0$. So we have shown that any element in
the coset of $x$ maps to $\delta_0$, modulo elements of filtration $2$. However, the
image of $V_F$ is contained in the homotopy derivations, by Lemma \ref{derivation}. By
Proposition \ref{coordfree}, there are no homotopy derivations of filtration higher than
one. Therefore the proposition is proved.
\end{proof}

\begin{cor}\label{firstalg}
For any regular triple $(R, T, F)$, there is a realization
\[
\vartheta_0\: D_{F/T} \to \Sigma I/I^2
\]
of $\delta_0\: F=T/I\to \Sigma I/I^2$ as a strict derivation over $T$.
\end{cor}

\begin{proof}
We have proved this for $R=T$. For the general case, Proposition \ref{strictbs}, applied
to the regular pair of $R$-algebras $(R, L)$, guarantees the existence of a strict
realization $\vartheta_0'\: D_L\to\Sigma J/J^2$ of $\delta'_0\: L\to \Sigma J/J^2$ from
Section \ref{toptowers}. Smashing it with the identity on $T$ produces an $F$-bilinear
map $T\smash D_L \to \Sigma I/I^2$. But $D_{F/T}$, as the fibre of $\mu\: F\smash_T F\to
F$, is equivalent to $T\smash D_L$, and hence we are done.
\end{proof}

\begin{cor}\label{firstinf}
The map $\rho\: T/I^2 \to T/I$ from the tower \eqref{quotiadic} can be realized as a map
of $R$-algebras under $T$.
\end{cor}

\begin{proof}
Form the singular extension associated to $\vartheta_0$.
\end{proof}

\begin{rem}
Note that $T/I^2$ is equivalent as an $R$-algebra to the smash product of $T$ with the
singular extension on the derivation $\vartheta_0'$ from the proof of Corollary
\ref{firstalg}.
\end{rem}

As a consequence of Corollary \ref{firstinf}, we obtain a cofibre sequence of
$T/I^2$-bimo\-dules of the form
\begin{equation}\label{ti2seq}
I/I^2 \lra T/I^2 \lra F \xra{\theta_0} \Sigma I/I^2.
\end{equation}
A chosen element $\phi\in(D/D^2)^\vee$ determines a map
\[
\alpha_*\: \mathscr D^*_{T/I^2}(F, \Sigma I/I^2) \cong \mathscr D^*_{T/I^2}(F,
F)\otimes_{F^*} D/D^2 \xra{1\otimes\phi}\mathscr D^*_{T/I^2}(F, F).
\]
This induces an $F^*$-linear map
\begin{equation}\label{hom2}
(D/D^2)^\vee \lra F^*_{T/I^2}(F), \quad \phi\mapsto \alpha_*(\theta_0).
\end{equation}
The composition with the natural map $F^*_{T/I^2}(F)\to F^*_R(F)$ corresponds to the
inclusion of the homotopy derivations $\DDer_R^*(F, F)$. In particular, \eqref{hom2} is
injective. Let $\wh\Ten_{F^*}(M^*)$ denote the completed tensor algebra on an
$F^*$-module $M^*$.

\begin{prop}\label{firstend}
The homomorphism \eqref{hom2} extends to an algebra isomorphism
\begin{equation}\label{algebraiso}
\wh\Ten_{F^*}^*((D/D^2)^\vee) \cong F^*_{T/I^2}(F).
\end{equation}
\end{prop}

\begin{proof}
Consider the universal coefficient spectral sequence
\[
E_2^{*,*} = \Ext^{*,*}_{T^*/I^2}(F^*, F^*) \Longrightarrow
\mathscr D^*_{T/I^2}(F, F).
\]
By Lemma \ref{ext2} below, the $E_2$-term is isomorphic as an algebra to the tensor
algebra $\Ten^*_{F^*}((D/D^2)^\vee)$. By \cite{sw}*{Prop.\@ 6.5} and the remarks above,
the image of \eqref{hom2} is represented by the first $\Ext$-group. Hence this consists
of permanent cycles. As the spectral sequence is multiplicative, this forces it to
collapse and hence to converge strongly. It follows that \eqref{algebraiso} induces an
isomorphism on the associated graded modules. As the spectral sequence converges
strongly, \eqref{algebraiso} itself is an isomorphism \cite{boardman}*{Thm.\@ 2.6}.
\end{proof}

\begin{cor}\label{surj}
The restriction map $F^*_{T/I^2}(F) \to F^*_T(F)$ is surjective.
\end{cor}

\begin{proof}
It is a map of algebras and hence given by the projection
\[
\wh\Ten^*_{F^*}((D/D^2)^\vee) \lra \Lambda^*_{F^*}((D/D^2)^\vee)
\]
(compare Remark \ref{remt}).
\end{proof}

\begin{lem}\label{ext2}
There is an isomorphism of algebras
\[
\Ext^{*,*}_{R^*/I^2}(F^*, F^*) \cong \Ten^*_{F^*}((D/D^2)^\vee).
\]
\end{lem}

\begin{proof}
It is easy to write down an $R^*/I^2$-free resolution of $F^*$. Namely, let $M^*$ be a
free $R^*/I^2$-module on generators $y_i$ of (cohomological) degree $-|x_i|$. There is an
$R^*/I^2$-free resolution of $F^*$ of the form
\begin{equation}\label{cc}
\cdots\xra{\delta_4} M^*\otimes M^* \otimes M^* \xra{\delta_3} M^*\otimes M^*
\xra{\delta_2} M^* \xra{\delta_1} R^*/I^2 \xra{\can} F^*,
\end{equation}
where $\otimes=\otimes_{R^*/I^2}$. The differential $\delta_s$ is defined as
\[
\delta_s(u_1\otimes\cdots\otimes u_s) = \bar\phi(u_1)\. u_2\otimes\cdots\otimes u_s,
\]
where $\bar\phi\: M^*\to R^*/I^2$ evaluates $y_i$ at $x_i$, so $\bar\phi(y_i) = x_i$. We
briefly explain where this complex comes from. Let $\RR^*=R^*[y_1, \ldots, y_n]$ be the
polynomial algebra on generators $y_i$ of degree as specified above. Let $\II$ be the
ideal $(y_1,\ldots,y_n)\lhd\RR^*$. View $R^*$ as an $\RR^*$-algebra via the homomorphism
of algebras $\phi\:\RR^*\to R^*$ defined by $\phi(y_i)=x_i$. Let $\FF^*$ and $\GG^*$ be
the $\RR^*$-algebras $\RR^*/\II$ and $\RR^*/\II^2$ respectively. Note that as an
$R^*$-module, $\FF^*$ is just $R^*$. As $\phi(\II)\subset I$,  $\phi$ induces compatible
homomorphisms $\FF^*\to F^*$ and $\GG^*\to G^*=R^*/I^2$, which in turn induce maps
\begin{equation}\label{dgamaps}
R^*\otimes^{\mathbf L}_{\RR^*} \FF^* \to F^*, \quad R^*\otimes^{\mathbf L}_{\RR^*} \GG^*
\to G^*
\end{equation}
in the homotopy category of differential graded $R^*$-algebras \cite{ss}*{§5}. Here we
use the standard notation $\otimes^{\mathbf L}_{\RR^*}$ for the left derived functor of
$\otimes_{\RR^*}$. The first map of \eqref{dgamaps} is an equivalence, as
\[
R^*\otimes^{\mathbf L}_{\RR^*} \FF^* \simeq \K_{\RR^*}(y_1, \ldots, y_n) \otimes_{\RR^*}
\FF^* \simeq \K_{R^*}(x_1, \ldots, x_n) \simeq F^*,
\]
where $\K$ denote Koszul complexes. Because $\II/\II^2$ is isomorphic to a sum of
suspensions of $\FF^*$, the natural map $R^*\otimes^{\mathbf L}_{\RR^*} \II/\II^2 \to
I/I^2$ is an equivalence as well. The diagram of cofibre sequences
\[
\xymatrix{R^*\otimes^{\mathbf L}_{\RR^*} \II/\II^2 \ar[r]\ar[d] & R^*\otimes^{\mathbf
L}_{\RR^*} \GG^* \ar[r]\ar[d] & R^*\otimes^{\mathbf L}_{\RR^*} \FF^* \ar[r]\ar[d] &
R^*\otimes^{\mathbf L}_{\RR^*} \Sigma  \II/\II^2 \ar[d]
\\
I/I^2 \ar[r] & G^* \ar[r] & F^* \ar[r] &  \Sigma I/I^2,}
\]
where all the maps are the natural ones, and the $5$-Lemma imply now that the second map
of \eqref{dgamaps} is also an equivalence.

Using the equivalences \eqref{dgamaps}, we derive an equivalence of homotopy differential
graded algebras
\begin{equation}\label{dgamap}
G^*\otimes^{\mathbf L}_{\GG^*}\FF^* \simeq (R^*\otimes^{\mathbf
L}_{\RR^*}\GG^*)\otimes^{\mathbf L}_{\GG^*}\FF^* \simeq R^*\otimes^{\mathbf
L}_{\RR^*}\FF^*  \simeq F^*.
\end{equation}
This allows us to realize $F^*$ as the bar construction $B^{R^*}(G^*, \GG^*, \FF^*)$
(this is analogous to \cite{ekmm}*{Prop.\@ IX.2.3}). The latter is homotopy equivalent to
the realization of the reduced bar resolution, which is of the form (with
$\otimes=\otimes_{R^*}$)
\begin{equation*}%\label{reducedbar}
\xymatrix{ \cdots \ar@<-1.5ex>[r]\ar@<-0.5ex>[r]\ar@<0.5ex>[r]\ar@<1.5ex>[r] & G^*\otimes
\II/\II^2 \otimes \II/\II^2 \otimes \FF^* \ar@<-1ex>[r]\ar[r]\ar@<1ex>[r] & G^*\otimes
\II/\II^2\otimes \FF^* \ar@<-0.5ex>[r]\ar@<0.5ex>[r] & G^*\otimes \FF^* .}
\end{equation*}
The chain complex associated to this simplicial resolution is precisely the chain complex
\eqref{cc} from before.

Mapping \eqref{cc} into $F^*$ kills all the differentials, so the statement is true
additively. One way to account for the multiplicative structure is by using a duality
argument. One checks that the resolution \eqref{cc} is compatible with the standard
comultiplication on $\Ten^*_{R^*}(M^*)$ and uses the fact that the multiplication on
$\Ext^{*,*}_{G^*}(F^*, F^*)$ is dual to the comultiplication on $\Tor_{*,*}^{G_*}(F_*,
F_*)$. Another possibility is to use the fact that the cosimplicial $R^*$-module obtained
by applying $\Hom^*_{G^*}(-, F^*)$ to the resolution $B_*^{R^*}(G^*, \GG^*, \FF^*)$
admits a composition product which is compatible with the one on $\Ext$-groups. This is
explained in more detail in the proof of a later result, Theorem \ref{maingen}, in an
analogous situation.
\end{proof}

\section{Topological Hochschild cohomology and ordinary cohomology}\label{thhcoh}

Let $F$ be a regular quotient algebra of $R$ and $G$ an $R$-algebra over $F$, with
structure map $\pi\: G\to F$. We can view $F$ as a $G$-bimodule via pullback along $\pi$
and consider $\THH_R^*(G, F)$. Forgetting the right $G$-action, we can also ask about the
left $G$-linear endomorphisms $F^*_G(F)$ of $F$. The aim of this section is to relate
these two invariants and to characterize a situation where we can determine both of
these. In a first part, we establish a connection between bimodule cohomology
$F^*_{F\smash F^\op}(-)$ and (left) module cohomology $F^*_F(-)$. In a second part, we
apply the results to compare $\THH_R^*(G, F)$ and $F^*_G(F)$ and to use this to compute
the two groups in a special situation.

\subsection{Bimodule cohomology}

Let $F$ be a regular quotient algebra of $R$ and let $\fenv=F\smash F^\op$ denote its
``enveloping algebra''. To compare $F^*_{\fenv}(M)$ and $F^*_F(M)$ for an $F$-bimodule
$M$, we are going to use the fact that $F^*_F(M)$ is a right comodule over the coalgebra
$F^*_R(F)$. The comultiplication on  $F^*_R(F)$ is induced by the product of $F$ in the
usual way. We have to be a bit careful in defining it, however, as we don't insist on $F$
being homotopy commutative. For this reason, we give some details.

\begin{prop}
The graded endomorphisms $F^*_R(F)$ admit a natural augmented $F^*$-coalgebra structure.
\end{prop}

\begin{proof}
The first thing to notice is that a priori several $F^*$-actions on $F^*_R(F)$ have to be
distinguished. There are two natural right and two natural left $F^*$-actions, defined in
an obvious way. Fortunately, they all agree in the present situation. The reason is that
the four $R^*$-actions obtained via pulling back along the surjection $\eta_*\: R_*\to
F_*$ necessarily coincide, because $R$ is commutative. We define the comultiplication as
the composition
\[
\Delta\: F^*_R(F) \xra{\mu^*} F^*_R(F\smash F) \cong F^*_R(F)\otimes_{F^*} F^*_R(F)
\]
of the map induced by the product $\mu$ on $F$ and a version of the K\"unneth isomorphism
\cite{sw}*{Prop.\@ 2.6} valid for non-commutative ring spectra. The counit is induced by
the unit $\eta$, by means of $\eta^*\: F^*_R(F)\to F^*_R(R)\cong F^*$. The augmentation
is defined as $\nu^*\: F^*\cong F^*_F(F) \to F^*_R(F)$, via ``restriction along $\eta$''.
It is a coalgebra map because $\Delta(\id_F) = \id_F \otimes \id_F$.
\end{proof}

Let $U\: \dfenv \to \mathscr{D}_F$ denote the ``restriction functor'' along the
$R$-algebra map $1\smash\eta\: F\cong F\smash R\to F\smash F$. Let
$\Cohom^*_{F^*_R(F)}(-,-)$ denote right $F^*_R(F)$-colinear maps and
$P(-)=\Cohom^*_{F^*_R(F)}(F^*, -)$ the functor of primitives.

\begin{prop}\label{coaction+}
The functor $F^*_F(U(-))$ takes values in right $F^*_R(F)$-comodules. There is a natural
transformation
\begin{equation}\label{nattrf}
\psi\: F^*_{\fenv}(-) \to P(F^*_F(U(-))).
\end{equation}
\end{prop}

\begin{proof}
For the first part, we define a natural map for left $F$-modules $X$
\begin{equation}\label{nattrans}
F^*_F(X)\otimes_{F^*} F^*_R(F) \lra F^*_F(X\smash F)
\end{equation}
precisely analogous to the K\"unneth transformation cited above, by sending an element
$f\otimes\phi$ to the left $F$-linear map
\[
X\smash F \xra{f\smash\phi} F\smash F \xra{\mu} F.
\]
As \eqref{nattrans} is an isomorphism on the suspensions of $F$, it is a natural
equivalence, by the usual comparison argument for cohomology theories. If $M$ is now an
$F$-bimodule, we compose the map $F^*_F(M) \to F^*_F(M\smash F)$ induced by the right
action of $F$ on $M$ with the inverse of \eqref{nattrans} for $X=U(M)$. It is
straightforward to check that the obtained natural map
\[
F^*_F(U(M)) \lra F^*_F(U(M)) \otimes_{F^*} F^*_R(F)
\]
defines a coaction of $F^*_R(F)$ on $F_F^*(U(M))$.

To define a natural transformation of the stated form, we associate to an element $f\in
F^*_{\fenv}(M)=[M, F]^*_{\fenv}$ the induced homomorphism of $F^*_R(F)$-comodules
\[
F^*\cong F^*_F(U(F)) \xra{F^*_F(U(f))} F^*_F(U(M)).
\]
This yields a natural map of the required form.
\end{proof}

Let $C^*$ be a graded coalgebra over a graded commutative ring $\Lambda^*$. Recall that a
graded right comodule $M^*$ over $C^*$ is called relatively injective if the functor
$\Cohom^*_{C^*}(-, M^*)$ takes split short exact sequences of $\Lambda^*$-modules to
short exact sequences (of $\Lambda^*$-modules). It is well-known that $M^*$ is relatively
injective if and only if $M^*$ is isomorphic to a retract of an extended $C^*$-comodule,
\ie\ a comodule of the form $X^*\otimes_{\Lambda^*} C^*$ for some $\Lambda^*$-module
$X^*$ (see \eg\ \cite{hoveyart}*{Lemma 3.1.2}). Recall that we may compute the derived
functors $\Coext^{*,*}_{C^*}(M^*, N^*)$ using relatively injective resolutions of $M^*$
if $M^*$ is $\Lambda^*$-projective \cite{ravenel}*{Lemmas A1.1.6(b), A1.2.9(b)}.

\begin{prop}
There are Adams-type spectral sequences
\[
\Coext^{*,*}_{F^*_R(F)}(F^*, F^*_F(X)) \ \Longrightarrow\ F^*_{\fenv}(X),
\]
whose edge homomorphism coincides with $\psi_X$ from \eqref{nattrf}.
\end{prop}

\begin{proof}
For given $X$, form the tower of $F$-bimodule spectra
\[
\xymatrix@C=0.2cm{ X=X_0 \ar[rr] && X_1 \ar[rr]\ar[ld]|-\circ && X_2\ar[ld]|-\circ
\ar[rr] &&  \cdots,
\\
& X_0\smash F \ar[lu] && X_1 \smash F \ar[lu] && \cdots\ar[lu] }
\]
where $X_{s+1}$ is inductively defined as the cofibre of right scalar multiplication on
$X_{s}$, $X_{s}\smash F \to X_{s}$. Applying $F^*_F(U(-))$ gives a relative injective
resolution of $F^*_F(U(X))$ over $F_R^*(F)$. Furthermore, for $X$ of the form $X=Y\smash
F$, we have $F^*_{\fenv}(X)\cong F^*_F(Y)$, and this isomorphism is induced by $\psi_X$
from \eqref{nattrf}:
\[
F^*_{\fenv}(X) \to \Hom^*_{F^*_R(F)}(F^*, F^*_F(U(X))) \cong F^*_F(Y).
\]
Therefore the spectral sequence exists by standard arguments.
\end{proof}

\begin{cor}\label{thekey}
If $M$ is an $F$-bimodule for which $F^*_F(U(M))$ is a relatively injective
$F^*_R(F)$-comodule, we have a natural identification $F^*_{\fenv}(M)\cong
P(F^*_F(U(M)))$.
\end{cor}

\begin{rem}
In a completely analogous manner, we can compare $F^*_{\fenv}(-)$ with $F^*_{F^\op}(-)$.
We then use the fact that $F^*_{F^\op}(U'(M))$ has a natural left $F^*_R(F)$-coaction for
an $F$-bimodule $M$, where $U' \: \dfenv\to \mathscr{D}_{F^\op}$ is the obvious
``restriction functor''.
\end{rem}

\begin{rem}
There is a dual version for homology of the previous statements. The situation is
considerably easier in homology, and we contend ourselves with a brief sketch. Consider
the algebra $\fenv_*=\pi_*(F\smash F^\op)$. The map $\mu_*\: \fenv_*\to F_*$ induced by
the product map is an algebra map, as a consequence of the fact that it is
$\fenv_*$-linear. Hence $\fenv_*$ is an augmented $F_*$-algebra. Let $Q(M_*)$ denote the
indecomposables $Q(M_*)=F_*\otimes_{\fenv_*} M_*$ of a left $\fenv_*$-module $M_*$. The
homotopy groups $M_*$ of an $F$-bimodule $M$ admit an obvious left $\fenv_*$-action. The
analogue of the natural transformation \eqref{nattrf} is now just given by the K\"unneth
map
\begin{equation}\label{nattrfhom}
Q(M_*) = F_*\otimes_{\fenv_*} M_* \lra F_*^{\fenv}(M).
\end{equation}
The existence of K\"unneth spectral sequences \cite{ekmm}*{Theorem IV.4.1} implies that
\eqref{nattrfhom} is an isomorphism whenever $M_*$ is relatively projective over
$\fenv_*$.
\end{rem}

\subsection{Topological Hochschild cohomology of algebras over regular
quotients}

We recall the following fundamental fact:

\begin{prop}[\cite{ekmm}*{Prop.\@ III.4.1}]\label{extscal}
Let $\phi\: R\to R'$ be a map of $\SS$-algebras. By pullback along $\phi$, we obtain a
functor $\phi^*\:\mathscr{M}_{R'}\to\mathscr{M}_R$. It has a left adjoint
$\phi_*\:\mathscr{M}_{R}\to\mathscr{M}_{R'}$, given by $\phi_*(M)= R'\smash_R M$. The
adjunction passes to one on derived categories.
\end{prop}

Recall the homotopy category $\hara$ of $R$-algebras over a fixed $R$-algebra $A$ from
Section \ref{singext}.

\begin{prop}\label{relate}
Let $F$ be a regular quotient algebra of $R$.
\begin{enumerate}
\item The functor $F^*_{?}(F)\: \harf\to \Mod_{F^*}$ takes
values in right $F^*_R(F)$-comodules.% algebras augmented over $F^*$.
\item There is natural transformation $\theta\: \THH^*_R(-, F) \to P(F^*_{?}(F))$. It is an
isomorphism for those $R$-algebras $G$ over $F$ for which $F^*_G(F)$ is relatively
injective over $F^*_R(F)$.
\end{enumerate}
\end{prop}

\begin{proof}
Let $G$ be an $R$-algebra over $F$ with structure map $\pi\: G\to F$.

(i) The adjunction of Proposition \ref{extscal} obtained from $\pi$ allows us to identify
$F^*_G(F)$ with $F^*_F(F\smash_G F)$. Now $F\smash_G F$ is an $F$-bimodule, and hence
$F^*_F(F\smash_G F)$ admits a right $F^*_R(F)$-coaction by Proposition \ref{coaction+}.

(ii) The algebra map $\pi\smash\pi\: G^{\text{e}}\to \fenv$ determines an adjunction
between $G^{\text{e}}$- and $\fenv$-modules, from which we deduce an isomorphism
\[
\THH_R^*(G, F)  = F^*_{G^{\text{e}}}(G) \cong F^*_{\fenv}(F\smash_G G\smash_G F) \cong
F^*_{\fenv}(F\smash_G F).
\]
We recall our standing convention that all functors are taken in the derived sense when
we work in homotopy categories. So we implicitly replace the $G$-bimodule $G$ cofibrantly
above. Now define $\theta_G$ as the composition
\[
\THH_R^*(G, F) \cong F^*_{\fenv}(F\smash_G F) \xra{\psi_{F\smash_G F}} P(F^*_F(F\smash_G
F)) \cong P(F^*_G(F)),
\]
where $\psi_{F\smash_G F}$ is the natural map from \eqref{nattrf}. The last statement is
a consequence of Corollary \ref{thekey}.
\end{proof}

To obtain information about $\THH_R^*(G, F)$, we are thus naturally led to consider
$F^*_G(F)$. We will prove the following fact together with Theorem \ref{maingen} below.

\begin{prop}\label{surjective}
The natural algebra map $p\: F^*_G(F) \to F^*_R(F)$ is surjective.
\end{prop}

We define the augmentation ideal $I(F^*_G(F))\subseteq F^*_G(F)$ to be the kernel of $p$:
\[
I(F^*_G(F)) = \ker(p\: F^*_G(F) \to F^*_R(F)).
\]
The map $p$ is $F^*_R(F)$-colinear, if we regard the coalgebra $F^*_R(F)$ as a right
comodule over itself. This is clear by construction of the right coaction of $F^*_R(F)$
on $F^*_G(F)$ and the comultiplication on $F^*_R(F)$. By left exactness of the primitives
functor, we deduce a short exact sequence of $F^*$-modules
\begin{equation}\label{ses1}
0 \lra P(I(F^*_G(F))) \lra P(F^*_G(F)) \lra F^* \lra 0.
\end{equation}
On the other hand, we have a short exact sequence
\begin{equation}\label{ses2}
0 \lra  \Der^{*-1}_R(G, F) \lra \THH^*_R(G, F) \lra  F^* \lra 0
\end{equation}
involving $\THH_R^*(G, F)$. It is obtained from the long exact sequence \eqref{les} for
$A=M=F$ and $B=G$, by noting that $\THH^*_R(G, F)\to F^*$ sends $\pi$ to $1$ and is
therefore surjective. The map $\theta_G\: \THH^*_R(G,F) \to P(F^*_G(F))$ from Proposition
\ref{relate} is compatible with the surjections of \eqref{ses1} and \eqref{ses2} by
construction. Hence it induces a map
\begin{equation}\label{augmtrf}
\wt\theta_G\: \Der^{*-1}_R(G, F) \lra P(I(F^*_G(F))).
\end{equation}
Recall the forgetful map $V^*_G\: \Der^*_R(G, F) \to \wt F^*_R(G)$ from Section
\ref{singext}.

\begin{prop}\label{propinattrf}
There is a natural map of $F^*_R(F)$-comodules
\begin{equation}\label{inattrf}
W^*_G\: I(F^*_G(F)) \lra \wt F^{*-1}_R(G),
\end{equation}
which makes the diagram
\[
\xymatrix@C=0.5cm{ \Der^{*-1}_R(G, F) \ar[rr]^-{\wt\theta_G} \ar[dr]_-{V^*_G} &&
P(I(F^*_G(F))) \ar[dl]^-{\ \ W^*_G|_{P(I(F^*_G(F)))}}
\\
& \wt F^{*-1}_R(G)}
\]
commutative.
\end{prop}

We will prove this together with Theorem \ref{maingen} as well.

\begin{prop}\label{algebraonthh}
The topological Hochschild cohomology $\THH_R^*(B, A)$ of an $R$-algebra $B$ over $A$ is
an $R^*$-algebra augmented over $A^*$, in a natural way.
\end{prop}

\begin{proof}
Let $\pi\: B\to A$ be the structure map of $B$. Applying Proposition \ref{extscal} to the
algebra map $1\smash\pi^\op\: B\smash B^\op \to B\smash A^\op$, we obtain:
\begin{equation}\label{first}
\THH^*_R(B, A) = A^*_{B^\text{e}}(B) \cong A^*_{B\smash A^\op}((B\smash
A^\op)\smash_{B^\text{e}} B).
\end{equation}
We claim that $(B\smash A^\op)\smash_{B^\text{e}} B\simeq A$ as $B\smash A^\op$-modules.
To prove this, we use the bar resolution $B_*^R(\tilde B,\tilde B,\tilde B)$ on a
cofibrant replacement $\tilde B$ of the $R$-module $B$ (see \cite{ekmm}*{IX.2}). This is
a simplicial strict $\tilde B^\text{e}$-module whose component of degree $q$ is $\tilde
B\smash (\tilde B)^{\smash q} \smash \tilde B$. The canonical map $\tilde B\smash \tilde
B \to \tilde B\smash_{\tilde B}\tilde B$ defines an augmentation and induces a derived
equivalence on applying geometric realization:
\[
B^R(\tilde B, \tilde B,\tilde B) = |B^R_*(\tilde B, \tilde B, \tilde B)| \simeq B\smash_B
B\simeq B.
\]
Now geometric realization commutes with smash products, and so we find
\begin{multline*}
(B\smash A^\op)  \smash_{B^\text{e}} B  \simeq (B\smash A^\op)\smash_{B^\text{e}}
|B^R_*(\tilde B, \tilde B, \tilde B)|
\\
\simeq |(\tilde B\smash \tilde A^\op)\smash_{\tilde B^\text{e}} B^R_*(\tilde B, \tilde B,
\tilde B)| \simeq |B^R_*(\tilde B, \tilde B, \tilde A)| \simeq B\smash_B A\simeq A,
\end{multline*}
where $\tilde A$ is a cofibrant replacement of the $R$-module $A$. Combining this
equivalence with \eqref{first}, we obtain
\[
\THH^*_R(B,A) \cong A^*_{B\smash A^\op}(A).
\]
Via this identification, the composition product on $A^*_{B\smash A^\op}(A)$ defines an
$R^*$-algebra structure on $\THH^*_R(B,A)$. The augmentation is induced by the natural
map of algebras $A^*_{B\smash A^\op}(A) \to A^*_{A^\op}(A) \cong A^*$, obtained by
pulling back along $\eta_B\smash 1\:A^\op \cong R\smash A^\op \to B\smash A^\op$.
\end{proof}

Here is the main result of this section:

\begin{thm}\label{maingen}
Let $F$ be a regular quotient algebra of $R$ and let $G$ be an $R$-algebra over $F$.
Assume that $W^*_G\: I(F^*_G(F)) \to \wt F^{*-1}_R(G)$ is surjective and that $F^*_R(G)$
is $F^*$-free.
\begin{enumerate}
\item There is an isomorphism of $F^*_R(F)$-comodules
\[
F^*_G(F) \cong \wh\Ten_{F^*}(\wt F_R^{*-1}(G))\otimes_{F^*} F^*_R(F).
\]
\item The forgetful map $V^*_G\:\Der^*_R(G, F)\to \wt F_R^*(G)$ is split surjective. Any chosen
section induces an isomorphism of algebras
\[
\wh\Ten_{F^*}(\wt F_R^{*-1}(G)) \cong  \THH_R^*(G, F).
\]
\end{enumerate}
\end{thm}

\begin{proof}
(i) In the following, we assume that $F$ and $G$ are both cofibrant as $R$-modules,
replacing them cofibrantly if necessary. Consider the two-sided bar resolution $B_*^R(G,
G, G)$ of $G$. It is of the form
\begin{equation}\label{bar}
\xymatrix{ \cdots  \ar@<-1ex>[r]\ar[r]\ar@<1ex>[r] &  G\smash G\smash G
\ar@<-0.5ex>[r]\ar@<0.5ex>[r] & G\smash G  \ar[r] & G.}
\end{equation}
On applying $F_{G}(-\smash_{G} F, F)$, we obtain a cosimplicial $R$-module equivalent to
one of the form
\begin{equation}\label{cosimplicial}
\xymatrix{F_{G}(F, F) \ar[r] & F_R(F, F) \ar@<-0.5ex>[r]\ar@<0.5ex>[r] & F_R(G\smash F,
F) \ar@<-1ex>[r]\ar[r]\ar@<1ex>[r] &\cdots}.
\end{equation}
It provides a resolution of $F_{G}(F, F)$ and leads to a (conditionally convergent)
Bous\-field--Kan spectral sequence
\begin{equation}\label{firstss}
E_r^{*,*} \Longrightarrow F^*_{G}(F).
\end{equation}
As we started with a resolution of $G$ as a bimodule over itself, the spectral sequence
is in fact one of $F^*_R(F)$-comodules. The unrolled exact couple underlying the spectral
sequence is of the form
\[
\xymatrix@C=0.2cm@R=0.5cm{ F^*_G(F) \ar[rd] && I(F^*_G(F)) \ar[ll]\ar[rd] && \cdots.
\ar[ll]
\\
& F^*_R(F) \ar[ru]|-\circ && F^{*-1}_R(G\smash F) \ar[ru]|-\circ }
\]
As the natural augmentations $F^*_G(F)\to F^*$ and $F^*_R(G^{\smash p}\smash F) \to F^*$
for $p\geq 0$ are all split and compatible, we deduce an exact couple involving reduced
groups
\begin{equation}\label{exactcouple}
\begin{array}{c}
\xymatrix@C=0.2cm@R=0.5cm{\wt F^*_G(F) \ar[rd]_-{j_0} && I(F^*_G(F))
\ar[ll]\ar[rd]_-{j_1} && \cdots.\ar[ll]
\\
&\wt F^*_R(F) \ar[ru]|-\circ &&\wt F^{*-1}_R(G\smash F) \ar[ru]|-\circ }
\end{array}
\end{equation}
We obtain a map $W^*_G\: I(F^*_G(F)) \to \wt F^{*-1}_R(G)$ as claimed in Proposition
\ref{propinattrf} by taking the composition
\[
W^*_G\: I(F^*_G(F))\xra{j_1} \wt F^{*-1}_R(G\smash F) \xra{(1\smash \eta_F)^*} \wt
F^{*-1}_R(G).
\]
The $E_1$-term of the spectral sequence is the cochain complex which is associated to the
cosimplicial $R^*$-module given by the homotopy groups of \eqref{cosimplicial}. As
$F^*_R(G)$ is free over $F^*$ by hypothesis, there are K\"unneth isomorphisms
\[
E_1^{p,*} = F^*_R((G)^{\smash p} \smash F) \cong F^*_R(G)^{\otimes_{F^*}(p)}\otimes_{F^*}
F^*_R(F).
\]
The components of the normalized cochain complex are given by
\[
N^{p,*}= \ker s^0\cap \cdots \cap \ker s^p \subseteq E_1^{p,*},
\]
where the $s^i$ are the codegeneracy maps. These are induced by
\[
G^{\smash(i+1)} \smash \eta_G \smash G^{\smash(p-i)} \smash F\: G^{\smash(p+1)} \smash F
\lra G^{\smash(p+2)} \smash F.
\]
It follows that the $s^i$ are given by
\[
a_1\otimes\cdots \otimes a_p\otimes x \mapsto \eta^*(a_i) a_1\otimes \cdots \otimes
\wh{a_i}\otimes\cdots\otimes a_p\otimes x
\]
where $\eta^*=F^*_R(\eta_G)\: F^*_R(G)\to F^*_R(R)\cong F^*$ is the augmentation and the
factor under $\wh{\ \ }\,$ is omitted. Hence we have
\[
N^{p,*} = \wt F^*_R(G)^{\otimes_{F^*} (p)} \otimes_{F^*} F^*_R(F).
\]
Now we use the fact that the spectral sequence has a multiplicative structure, which
corresponds to the composition pairing on the target. It is induced by
\[
F_R(G^{\smash p}\smash  F, F) \smash F_R(G^{\smash q}\smash F, F) \lra
F_R(G^{\smash(p+q)}\smash F, F)
\]
which sends $\alpha\smash\beta$ to the composition
\[
G^{\smash(p+q)} \smash F \xra{G^{\smash p} \smash \beta} G^{\smash p} \smash
F\xra{\alpha} F.
\]
The proof that this gives rise to a product structure on the spectral sequence is
analogous to the one given in \cite{bk} for the homotopy spectral sequence of a space
with coefficients in a ring. In particular, we find that our spectral sequence has an
action of $E_1^{0,*} \cong F^*_R(F)$. It is not difficult to see that the (reduced)
$E_1$-term is given as an algebra as
\[
\Ten^{*,*}_{F^*}(\wt F_R^*(G)) \otimes_{F^*} F^*_R(F).
\]
As the augmentation $F^*_G(F)\to F^*$ is surjective, the unit element $1\in E_1^{0,*}$
must be a permanent cycle. By $F^*_R(F)$-linearity, this implies that the whole of
$E_1^{0,*}$ consists of permanent cycles. As a consequence, we find that $p\: F^*_G(F)
\to F^*_R(F)$ is surjective, as claimed in Proposition \ref{surjective}. (This argument
is independent of the hypotheses made on $F^*_R(G)$ and $W^*_B$.) Now $W_G^*$ is
surjective by assumption, and hence $\wt F^*_R(G)$ consists of permanent cycles. The
multiplicative structure now forces the whole spectral sequence to collapse, and so it
converges strongly. This means that the canonical map
\[
F^*_G(F) \lra \lim_s F^*_G(F)/F^s(F^*_G(F))
\]
is an isomorphism, where
\[
\cdots \subseteq F^{s+1}(F^*_G(F)) \subseteq F^s(F^*_G(F)) \subseteq \cdots  \subseteq
F^0(F^*_G(F)) = F^*_G(F)
\]
is the spectral sequence filtration. As the filtration subquotients
\[
F^s(\cdots)/F^{s+1}(\cdots) \cong E_\infty^{s,*} \cong \wt F^*_R(G)^{\otimes_{F^*}(s)}
\otimes_{F^*} F^*_R(F)
\]
are all relatively injective, there are no non-trivial extensions. Furthermore, forming
extended comodules commutes with inverse limits, by definition of inverse limits of
comodules (see \cite{hoveyart}*{Prop.\@ 1.2.2}). This finishes the proof of part (i).

(ii) Part (i) and Proposition \ref{relate} imply that
\begin{equation}\label{moduleiso}
\THH_R^*(G, F) \cong \wh\Ten_{F^*}(\wt F^{*-1}_R(G))
\end{equation}
as $F^*$-modules. To obtain more information, we apply $F_{G\smash G^\op}(-, F)$ to the
bar resolution \eqref{bar}. This yields a cosimplicial $R$-module of the form
\begin{equation}\label{cosimplicial2}
\xymatrix{\THH_R(G, F) \ar[r] & R \ar@<-0.5ex>[r]\ar@<0.5ex>[r] & F_R(G, F)
\ar@<-1ex>[r]\ar[r]\ar@<1ex>[r] & F_R(G\smash G, F)
\ar@<-1.5ex>[r]\ar@<-0.5ex>[r]\ar@<0.5ex>[r]\ar@<1.5ex>[r] & \cdots}.
\end{equation}
Let us consider the associated Bousfield--Kan spectral sequence
\begin{equation}\label{secondss}
\bar E_r^{*,*} \Longrightarrow \THH_R^*(G, F).
\end{equation}
Just as the other spectral sequence \eqref{firstss}, it carries a multiplicative
structure, which is compatible with the algebra structure on the target from Proposition
\ref{algebraonthh}. Using similar arguments as above, we find that the reduced $E_1$-term
is isomorphic as an algebra to $\Ten_{F^*}(\wt F^*_R(G))$.

Now we assert that there is a morphism of spectral sequences $\omega\: \bar E_r^{*,*} \to
E_r^{*,*}$ compatible with $\theta_G\: \THH_R^*(G, F)\to F^*_G(F)$ on the targets. To
construct $\omega$, we note that the resolution giving rise to \eqref{secondss} can also
be obtained by applying $F_{\fenv}(F\smash_G - \smash_G F, F)$ to \eqref{bar}. Similarly,
the one inducing \eqref{firstss} is equivalent to the one obtained by applying
$F_{\fenv}(F\smash_G - \smash_G F\smash F, F)$ to \eqref{bar}. The product map of $F$
induces a morphism $\omega$ of the required form. On reduced $E^1$-terms, $\omega$ is
injective, as it is given by
\[
1^{\otimes p} \otimes \nu^*\: \wt F^*_R(G)^{\otimes_{F^*}(p)} \lra \wt
F^*_R(G)^{\otimes_{F^*}(p)}\otimes_{F^*} F^*_R(F),
\]
where $\nu^*$ is the augmentation of the coalgebra $F^*_R(F)$. It follows that the
spectral sequence \eqref{secondss} collapses, along with \eqref{firstss}. In particular,
we obtain a surjection
\begin{equation}\label{surje}
\Der^{*-1}_R(G, F) \cong \ker\bigl(\THH_R^*(G, F)\to F^* \bigr) \lra \bar
E_\infty^{1,*-1} = \wt F^{*-1}_R(G).
\end{equation}
This map coincides with the forgetful map $V^*_{G}$, as we show next. This will prove the
first statement of (ii), as $F^*_R(G)$ is $F^*$-free by assumption and as $\THH^*_R(G,
F)$ and hence $\Der^*_R(G, F)$ are $F^*$-modules by \eqref{moduleiso}. It will also imply
Proposition \ref{propinattrf}.

By \cite{laz}*{Thm.\@ 4.1} the universal derivation $d_{G}\: G \to D_{G}$ is the unique
map rendering the diagram
\[
\xymatrix@C=1.5cm{G \ar[d]^-{1\smash\eta\smash 1} \ar@{..>}[r] & D_{G}\ar[d]^-{\can}
\\
G\smash G\smash G \ar[r]^-{\mu\smash 1- 1\smash\mu} & G\smash G}
\]
commutative. Now the first bit of the tower which underlies the spectral sequence
\eqref{secondss} is of the form
\begin{equation}\label{tower}
\begin{array}{c}
\xymatrix@C=0.3cm{\cdots \ar[rd] &&  D_{G}\ar[rd]^-\can \ar[ll]|-\circ && G.
\ar[ll]|-\circ
\\
&G \smash G\smash G \ar[rr]^-{\mu\smash 1-1\smash \mu} \ar[ru]^-{e_G} &&G\smash G
\ar[ru]_-\mu}
\end{array}
\end{equation}
The map labeled $e_G$ is uniquely determined by the requirement that it be a lift of
$\can\: D_G \to G\smash G$. Now $d'_G\: G\smash G\smash G \to D_G$ obtained from $d_G$ by
extending scalars is such a lift, hence $e_G=d'_G$. By definition of $V_G^*$, this shows
that \eqref{surje} is $V^*_G$.

To prove the last statement, we choose an arbitrary section to $V^*_G$ and consider the
induced map of algebras
\[
\wh\Ten_{F_*}(\wt F^{*-1}(G)) \lra \THH_R^*(G, F).
\]
It is an isomorphism on the associated graded and therefore an isomorphism itself, by
\cite{boardman}*{Thm.\@ 2.6}.
\end{proof}

\section{Higher order infinitesimal thickenings of $F$}\label{infthick}

We will now prove the following result, which easily implies Theorem \ref{theorem} as we
show below:

\begin{thm}\label{theoremplus}
Let $F$ be a regular quotient algebra of an even commutative $\SS$-algebra $R$. There is
a sequence of $R$-algebra spectra under $R$ of the form
\begin{equation*}
R\lra \cdots \lra R/I^{s+1}\lra R/I^s \lra \cdots \lra R/I=F,
\end{equation*}
whose sequence of coefficients is the canonical sequence of $R_*$-algebras
\[
R_* \lra \cdots \lra R_*/I^{s+1}\lra R_*/I^s \lra \cdots \lra R_*/I \cong F_*
\]
and for which the homomorphisms $W_{R/I^s}^*\: I(F^*_{R/I^s}(F)) \to \wt
F^{*-1}_R(R/I^s)$ are surjective for $s>1$.
\end{thm}

As $F^*_R(R/I^s)$ is $F^*$-free by Proposition \ref{ses}, Theorem \ref{maingen} then
implies:

\begin{cor}\label{thhregular}
The forgetful map $V^*_{R/I^s}\: \Der^*_R(R/I^s, F) \to \wt F^*_R(R/I^s)$ is split
surjective. Any splitting gives rise to an algebra isomorphism
\[
{\wh \Ten}^*_{F}(\widetilde{F}^{*-1}_R(R/I^s)) \cong \THH_R^*(R/I^s, R/I).
\]
Furthermore, there is an isomorphism of $F^*_R(F)$-comodules
\[
F^*_{R/I^s}(F) \cong \wh\Ten_{F^*}(\wt F_R^{*-1}(R/I^s))\otimes_{F^*} F^*_R(F).
\]
\end{cor}

\begin{proof}[Proof of Theorem \ref{theoremplus} $\Longrightarrow$ Theorem \ref{theorem}]
By definition of a regular pair, $F$ is of the form $T\smash L$, where $L=R/J$ is a
regular quotient algebra of $R$.  Let
\[
R \lra \cdots \lra R/J^{s+1} \lra R/J^s \lra \cdots \lra R/J = L
\]
be a sequence of $R$-algebras with the properties specified in Theorem \ref{theoremplus}.
By applying $T\smash-$, we obtain a sequence of $R$-algebras under $T$ satisfying the
required conditions (see the proof of \cite{sw}*{Thm.\@ 6.6} for details).
\end{proof}

\begin{proof}[Proof of Theorem \ref{theoremplus}]
We construct the sequence inductively. The start of the induction is provided by
Corollary \ref{firstinf}. Assuming that the sequence is constructed up to $R/I^s$, we
show now that $W_{R/I^s}^*\: I(F^*_{R/I^s}(F)) \to \wt F^{*-1}_R(R/I^s)$ is surjective.
Consider the map $\pi_s\: R/I^s \to F$ given by composing all the maps of the sequence
constructed so far and regard it as a map in the derived category of $R/I^s$-bimodules.
Let
\begin{equation}\label{cof1}
I/I^s \xra{\iota_s} R/I^s \xra{\pi_s} F \xra{\tau_s} \Sigma I/I^s
\end{equation}
be a cofibre sequence associated to $\pi_s$ and let
\begin{equation}\label{cof2}
\Sigma^{-1} R/I^s \lra D_{R/I^s}\lra R/I^s\smash R/I^s \xra{\mu_s} R/I^s
\end{equation}
be a cofibre sequence associated to the product $\mu_s$ of $R/I^s$. Smashing the latter
over $R/I^s$ with $F$ and $R/I^s$ from the right respectively gives rise to a diagram (we
write $\smash_s$ for $\smash_{R/I^s}$ and $D_s$ for $D_{R/I^s}$)
\[
\xymatrix{ \Sigma^{-1} F \ar[r] & D_s\smash_s F \ar[r] & R/I^s\smash F \ar[r] & F
\\
\Sigma^{-1} R/I^s \ar[r] & D_s \ar[r]\ar[u]^-{1\smash\pi_s} & R/I^s \smash R/I^s
\ar[r]\ar[u]^-{1\smash\pi_s} & R/I^s}
\]
By the $3\times 3$-Lemma \cite{hps}*{Lemma A.1.2}, we can complete it to one of the form
\begin{equation}\label{bigdiagram}\begin{array}{c}
\xymatrix{ \Sigma^{-1} X_s \ar[r] & D_s\smash_s \Sigma I/I^s \ar[r] & R/I^s \smash \Sigma
I/I^s \ar[r] & X_s
\\
\Sigma^{-1} F \ar[r]\ar[u] & D_s\smash_s F \ar[r]\ar[u]^-{1\smash\tau_s} & R/I^s\smash F
\ar[r]\ar[u]^-{1\smash\tau_s} & F\ar[u]
\\
\Sigma^{-1} R/I^s \ar[r]\ar[u] & D_s \ar[r]\ar[u]^-{1\smash\pi_s} & R/I^s \smash R/I^s
\ar[r]\ar[u]^-{1\smash\pi_s} & R/I^s\ar[u]
\\
\Sigma^{-2} X_s \ar[r]\ar[u] & D_s\smash_s I/I^s \ar[r]\ar[u]^-{1\smash\iota_s} &
R/I^s\smash I/I^s \ar[r] \ar[u]^-{1\smash\iota_s} & \Sigma^{-1} X_s, \ar[u] }
\end{array}
\end{equation}
for some $X_s$, such that all squares commute, except for the bottom left one, which
anti-commutes, and such that all the rows and columns are cofibre sequences. We apply
$F^*_s(-)=F^*_{R/I^s}(-)$ and obtain
\[
\def\objectstyle{\scriptstyle}
\xymatrix@R=0.7cm@C=0.7cm{ & \cdots\ar[d] & \cdots\ar[d] & \cdots\ar[d] & \cdots\ar[d]
\\
\cdots & \ar[l] F^{*+1}_s(X_s)\ar[d] & F^*_s(D_s\smash_s \Sigma I/I^s) \ar[l] \ar[d]&
F^{*-1}_R(I/I^s) \ar[l] \ar[d] & F^*_s(X_s) \ar[l]\ar[d] & \cdots \ar[l]
\\
\cdots & \ar[l] F^{*+1}_s(F) \ar[d]  & F^*_s(D_s\smash_s F)\ar[d] \ar[l] & F^*_R(F)
\ar[l]\ar[d] \ar@{}[dr]|*{(\ast)} & F^*_s(F)\ar[l]\ar[d]& \cdots \ar[l]
\\
\cdots & \ar[l] F^{*+1} \ar[d]& F^*_s(D_s) \ar[l] \ar[d]& F^*_R(R/I^s) \ar[l] \ar[d]& F^*
\ar[l]\ar[d]& \cdots \ar[l]
\\
\cdots & \ar[l] F^{*+2}_s(X_s) \ar[d]& F^*_s(D_s\smash_s I/I^s) \ar[l]\ar[d] &
F^{*}(I/I^s) \ar[l]\ar[d] & F^{*+1}_s(X_s) \ar[l]\ar[d] & \cdots \ar[l]
\\
& \cdots & \cdots & \cdots & \cdots }
\]
The fact that square $(\ast)$ commutes implies that the vertical map $R/I^s\to F$ in the
right-most column of \eqref{bigdiagram} is the map $\pi_s$ from above, when viewed as a
left $R/I^s$-linear map. Therefore this column is equivalent to the cofibre sequence
\eqref{cof1} of left $R/I^s$-modules. Now observe that the map canonical $F^*_s(F) \to
F^*_R(F)$ is surjective, because it factors through the surjection $F^*_{R/I^2}(F) \to
F^*_R(F)$ (Corollary \ref{surj}). We claim that $F^{*+1}_s(X_s)\cong F^*_s(I/I^s) \to
F^{*}_R(I/I^s)$ is surjective, too. To show this, we ``invert'' the sequence $R/I^s \to
\cdots \to R/I^2 \to F$ (of $R/I^s$-bimodules) by repeated use of the octahedral axiom to
obtain a sequence
\[
I/I^s \lra I/I^{s-1} \lra \cdots \lra I/I^3 \lra I/I^2.
\]
We show inductively that the natural maps $F^*_s(I/I^t)\to F^*_R(I/I^t)$ are surjections
for $t\leq s$. This is the case for $t=2$, as $I/I^2$ is a wedge of suspension of $F$.
Assume that $F^*_s(I/I^t)\to F^*_R(I/I^t)$ is surjective. Note that the homotopy fibre of
$I/I^{t+1}\to I/I^t$ is equivalent to the one of $R/I^{t+1}\to R/I^t$, which is
$I^t/I^{t+1}$ by construction. Consider the exact sequences
\[
\xymatrix{ \cdots \ar[r] & F^*_s(I/I^t) \ar[r]\ar[d] & F^*_s(I/I^{t+1}) \ar[r] \ar[d] &
F^*_s(I^t/I^{t+1}) \ar[r] \ar[d] & F^{*+1}_s(I/I^t) \ar[r] \ar[d] & \cdots
\\
\cdots \ar[r] & F^*_R(I/I^t) \ar[r]& F^*_R(I/I^{t+1}) \ar[r]& F^*_R(I^t/I^{t+1}) \ar[r]
 & F^{*+1}_R(I/I^t) \ar[r] & \cdots}
\]
As $I^t/I^{t+1}$ is given as a wedge of suspensions of $F$, the inductive step follows
from the $5$-Lemma. What we have seen implies, together with Proposition \ref{ses}, that
we can extract from \eqref{bigdiagram} a diagram with exact rows and sequences of the
form
\[
\xymatrix@R=0.5cm{& 0\ar[d]
\\
& \wt F^{*-1}_R(R/I^s) \ar[d] & 0 \ar[d]
\\
0 & \ar[l] F^{*-1}_R(I/I^s) \ar[d] & F^*_s(X_s)\ar[d] \ar[l] & F^{*-1}(D_s\smash_s \Sigma
I/I^s) \ar[d]\ar[l] & 0 \ar[l]
\\
0 & \ar[l] \wt F^*_R(F)\ar[d] & \wt F^*_s(F) \ar[l]\ar[d] & F^{*-1}_s(D_s\smash_s
F)\ar[d]^-{(1\smash\pi_s)^*}\ar[l] & 0 \ar[l]
\\
& 0 & 0 & F^{*-1}_s(D_s)}
\]
The snake lemma shows that the upper vertical map in the right column is injective, the
lower one, $(1\smash\pi_s)^*$, is surjective and that $F^*_s(D_s)\cong \wt F^*_R(R/I^s)$.
We show now that we can factorize $W_{R/I^s}^*\:I(F^*_s(F))\to\wt F^{*-1}_R(R/I^s)$ as
\begin{equation}\label{desired}
I(F^*_s(F)) \cong F^{*-1}_s(D_s\smash_s F) \xra{(1\smash\pi_s)^*} F^{*-1}_s(D_s)
\xra{\psi} \wt F^{*-1}_R(R/I^s),
\end{equation}
where $\psi$ is an isomorphism. This will imply surjectivity of $W_{R/I^s}^*$.

Recall the definition of $W_{R/I^s}^*$, based on the exact couple \eqref{exactcouple}
(for $G=R/I^s$). This exact couple can be obtained by applying $F^*_s(-)$ to a diagram of
the form
\[
\xymatrix@C=0.2cm@R=1cm{ F \ar[rr]|-\circ  && D_s\smash_s F \ar[ld]\ar[rr]|-\circ &&
\cdots\ar[ld]
\\
& R/I^s\smash F \ar[lu] && R/I^s\smash R/I^s\smash F\ar[lu]^-{\phi}   }
\]
and taking reduced cohomology where appropriate. Arguing as after for \eqref{tower}, we
find that $\phi$ must be the map obtained from the universal derivation $d_s\: R/I^s \to
D_s$ by first extending scalars and then applying $-\smash_s F$. Therefore $W^*_{R/I^s}$
is given as the composition
\begin{align*}
I(F^*_s(F))\cong & F^{*-1}_s(D_s\smash_s F) \xra{\phi^*} \wt F^{*-1}_s(R/I^s\smash R/I^s
\smash F) \xra{(1\smash 1\smash \eta_F)^*}
\\
& \wt F^{*-1}_s(R/I^s \smash R/I^s) \cong \wt F^{*-1}_R(R/I^s).
\end{align*}
Denoting by $\phi'\: R/I^s\smash R/I^s \to D_s$ the map induced from $d_s$ by extending
scalars to $R/I^s$ on the left, we have a commutative diagram
\[
\xymatrix{R/I^s \smash R/I^s \ar[r]^-{\phi'} \ar[d]^-{1\smash 1\smash \eta_F} & D_s
\ar[d]^-{1\smash \pi_s}
\\
R/I^s \smash R/I^s \smash F \ar[r]^-\phi & D_s \smash_s F.}
\]
It allows us to factor $W^*_{R/I^s}$ as
\begin{align*}
I(F^*_s(F)) &\cong F^{*-1}_s(D_s\smash_s F) \xra{(1\smash \pi_s)^*} \wt F^{*-1}_s(D_s)
\xra{(\phi')^*} \wt F^{*-1}_s(R/I^s \smash R/I^s)
\\
& \cong \wt F^{*-1}_R(R/I^s).
\end{align*}
It is straightforward to check that the map $D_s\to R/I^s\smash R/I^s$ from the fibre
sequence \eqref{cof2} induces an inverse to $(\phi')^*$. Therefore, we have found a
factorization of $W^*_{R/I^s}$ of the form \eqref{desired} and shown that it is
surjective.

We now construct $R/I^{s+1}$ as an $R$-algebra. Theorem \ref{maingen} implies that
\[
V_{R/I^s}^* \: \Der_R^*(R/I^s, \Sigma I^s/I^{s+1}) \to \wt{(\Sigma I^s/I^{s+1})}^*(R/I^s)
\]
is surjective, as $I^s/I^{s+1}$ is a wedge of suspensions of $F$. Observe that the map
$\theta^s\: R/I^s \to \Sigma I^s/I^{s+1}$ from the tower \eqref{quotiadic} is contained
in the reduced cohomology group
\[
(\wt{I^s/I^{s+1}})^*_R(R/I^s) = \ker\bigl((I^s/I^{s+1})^*_R(R/I^s)\xra{\eta_{R/I^s}^*}
(I^s/I^{s+1})^*_R(R) \cong I^s/I^{s+1}\bigr)
%
%\cong I^s/I^{s+1}\otimes_{F^*} \wt F_R^*(R/I^s).
\]
This is clear, as $R/I^{s+1}\xra{\rho} R/I^s\xra{\theta^s} \Sigma I^s/I^{s+1}$ is zero
and as the unit $\eta_{R/I^s}$ factors through $\rho$. Thus we can realize $\theta_s$ as
a strict derivation and therefore construct $R/I^{s+1}$ as the associated singular
extension.
\end{proof}

\begin{rem}
If $R$ is the Eilenberg-MacLane spectrum $H\Lambda$ of a commutative ring $\Lambda$,
Corollary \ref{thhregular} states that
\begin{equation}\label{algebraic}
\Ext^*_{\Lambda/I^s}(\Lambda/I, \Lambda/I) \cong
T^*_{\Lambda/I}(\wt\Ext_{\Lambda}^{*-1}(\Lambda/I^s, \Lambda/I)) \otimes_\Lambda
\Ext^*_{\Lambda}(\Lambda/I, \Lambda/I),
\end{equation}
by \cite{ekmm}*{IV.2}. This is a well-known fact in commutative algebra. It is usually
considered in the case where $\Lambda$ is a regular local ring with maximal ideal $I$.
Equation \eqref{algebraic} implies that $\Lambda/I^s$ is a Golod ring under this
assumption. See \cite{avramov} for details.
\end{rem}

\end{document}